\documentclass[final]{elsarticle}
\usepackage{latexsym, graphicx, epsfig, amsmath, amsfonts,amssymb,subfigure,algorithm,algorithmicx,booktabs,bm,color}

\def\x{\xi}
\def\Pb{\textrm{Prob}}
\def\mb{\mathbf}

\def\z{\zeta}

\def\i{\mathbf{i}}
\def\j{\mathbf{j}}

\def\t{\tilde}

\def\E{{\mathbb E}}
\def\R{{\mathbb R}}
\def\P{{\mathcal P}}

\def\O{{\Omega}}
\def\o{{\omega}}
\def\A{{\mathcal{A}}}

\def\_#1{{\underline{#1}}}

\newproof{proof}{Proof}
\begin{document}
\begin{frontmatter}

\title{A unified framework for mesh refinement in random and physical space}

\author{Jing Li}

\address{Pacific Northwest National Laboratory. Email:jing.li@pnnl.gov.}

\author{Panos Stinis}

\address{Pacific Northwest National Laboratory. Email:panagiotis.stinis@pnnl.gov.}

\begin{abstract}
In recent work we have shown how an accurate reduced model can be utilized to perform mesh refinement in random space. That work relied on the explicit knowledge of an accurate reduced model which is used to monitor the transfer of activity from the large to the small scales of the solution. Since this is not always available, we present in the current work a framework which shares the merits and basic idea of the previous approach but {\it does} not require an explicit knowledge of a reduced model. Moreover, the current framework can be applied for refinement in both random and physical space. In this manuscript we focus on the application to random space mesh refinement. We study examples of increasing difficulty (from ordinary to partial differential equations) which demonstrate the efficiency and versatility of our approach. We also provide some results from the application of the new framework to physical space mesh refinement.   
\end{abstract}

\end{frontmatter}

\section{Introduction} \label{sec:intro}
The prediction of the evolution of most physical systems is inevitably affected by uncertainties in the initial, boundary and/ or parameters entering the formulation of the governing equations. The two prevalent methods for quantifying the effect of this uncertainty on the evolution of the statistics are Monte Carlo (MC) and stochastic spectral method (see, e.g., \cite{XiuK_SISC02,XiuH_SISC05,Xiu10,BabuskaTZ_SINUM04,DebBO01,GhanemS91,FrauenfelderST_CMAME05,Schwab_IMA07,Babuska2007,MathelinH_NASA03,MathelinHZ_NA03,MatthiesK_CMAME05} and references therein). Both of these methods lead to efficient algorithms as long as the number of sources of uncertainty (uncertainty dimensionality) is relatively small. However, for most realistic systems, the dimensionality of uncertainty is rather large. This can lead to very serious obstacles in the application of both MC and stochastic spectral method. Development of sparse grid techniques \cite{Keese2003,XiuH_SISC05,nobile2008sparse,nobile2008anisotropic} as well as sparse polynomial bases \cite{Schwab_IMA07} for the Galerkin formulation have greatly alleviated this difficulty, however these techniques heavy rely on the regularity of the solution in the random space. At the same time, not all values in the range of the sources of uncertainty contribute equally to the determination of the solution statistics. This creates the hope that if one can allocate more computational resources around the most important values of the random variables, then one can still produce reasonably accurate predictions of the solution statistics. Multi-element generalized polynomial chaos method (ME-gPC) \cite{WanK_JCP05,WanK_SISC06} multi-element probabilistic collocation method (ME-PCM) \cite{Foo2008} and multi-resolution wavelet expansion \cite{MaitreKNG_JCP04,MaitreNGK_JCP04} provide the framework to obtain the statistics of the system when applying different computational resources to different area (element) in random space. In order to affect the concentration of resources in the more sensitive areas of the random space (for example, the discontinuities in the random space), one needs tools for the reliable location of these sensitive areas. In particular, if we consider a mesh discretization of the random space, then one needs to be able to decide when, where and in which direction to refine this mesh. \cite{WanK_JCP05} provides an approach to refine the random elements and allocate the computational resources adaptively as the system evolves. Adaptive sparse grid methods \cite{Ma2009,jakemanG2013,Zhang2013} also provided a dynamical way to obtain unconstructed collocation mesh to resolve the system. 

In \cite{Li2015} we presented a novel way of performing mesh refinement in random space. The approach was based on the use of a reduced model of the full order model, when such a reduced model is available in explicit form. The main idea was that successful mesh refinement depends on the accurate monitoring of the transfer of activity (e.g. mass, energy) from the large to the small scales of the solution. At the same time, this is exactly the role of a good reduced model, to reproduce faithfully this transfer of activity. Thus, if one has access to a good reduced model, then one can use it to perform the mesh refinement task. In \cite{Li2015} we provided theoretical results which established the soundness of our approach as well as numerical experiments which corroborated the theory. However, the whole framework relied on the explicit knowledge of the reduced model which is not always easy to obtain.

In the current work, we present a simple reformulation of the idea in \cite{Li2015} which retains the merits of the previous approach while at the same time eliminating the need for the explicit knowledge of the reduced model. In particular, we offer a way to compute the rate of transfer of activity from the large to the small scales by using appropriate quantities which can be computed {\it without} the explicit knowledge of the reduced model (see Section \ref{sec:mesh_refinement}). The main idea is that at every step and at every element where refinement may be needed, one can obtain an expansion of the solution (in the element) on some chosen basis. The coefficients of this expansion can be used to compute the rate of transfer of activity from the large to the small scales (again within the element) and thus, decide whether there exists a need for refinement. We present results for the uncertainty quantification of three different cases of increasing complexity (linear ODE with uncertain coefficient, Kraichnan-Orszag system and Kuramoto Sivashinsky equation with an uncertain bifurcation parameter). In all three cases, the mesh refinement algorithm is capable of detecting the sensitive areas of the random space with high accuracy. 

The adaptive mesh refinement framework described above can also be applied to situations where refinement is required in physical space. The only difference with the random space case is that in physical space when assigning values for the modes in the newly generated elements, the coupling between adjacent elements is required. In other words, the elements are not decoupled and the continuity (or more) of the solution across the boundary of adjacent elements should be guaranteed. In the current work we present numerical results for the mesh refinement around the location of shocks for the 1D Burgers equation while more elaborate applications will appear in a forthcoming publication.


\section{Representation of uncertainty} \label{sec:gPC}
Let $(\O,\A,\P)$ be a complete probability space, where $\O$ is the event space and $\P$ is the probability measure defined on $\A\in2^{\O}$, the $\sigma-$ algebra of subsets of $\O$ (these subsets are called events). 
Let $D$ be a subset of $\R^c(c\in\{1,2,3\})$ with boundary $\partial D$. Let $\mathcal{L}$ be an operator on $D$, and
consider the following stochastic differential equation:
\begin{equation}\label{sto_eq}
u_t(\mb{x},t;\o) = \mathcal{L}(\mb{x},t,\o;u), \qquad \mb{x}\in D.
\end{equation}
The operator $\mathcal{L}$ usually involves spatial derivatives and can be either linear or nonlinear. In addition, $\mathcal{L}$ may depend on $\o\in\O$. Appropriate initial conditions and boundary conditions sometimes involving random parameters relying on $\o$ are assumed so that the problem \eqref{sto_eq} is well-posed $\P-$a.e. We also assume that for $\P-$a.e. solution $u:=u(\mb{x},t;\o)$ is a function taking values in $\R$. 

The random dependence operator $\mathcal{L}$ must satisfy a few important properties. The most important one is the "finite dimensional noise assumption" \cite{Babuska2007,XiuH_SISC05}, that is the random input can be represented with a finite-dimensional probability space. More precisely, the random input can be represented by a finite set of random variables, e.g., $\bm{\x}(\o) = (\x_1(\o),\cdots,\x_d(\o))$ be a $d$-dimensional random vector for $\o\in\O$. Then by the Doob-Dynkin lemma \cite{Oksendal98}, the solution $u(\mb{x},t;\o)$ can be written as $u(\mb{x},t;\bm{\x}(\o))$. 

Any second-order stochastic process can be represented as a random variable at each spatial and temporal location. Applying the Karhunen-Lo${\grave{e}}$ve expansion \cite{GhanemS91} along with the "finite dimensional noise assumption", a second-order stochastic process can be characterized by a finite set of mutually independent random variables.  Thus without loss of generality, in this work we assume $\{\xi_i(\o)\}_{i=1}^N$ are independent random variables. Let $\Gamma:=\prod_{i=1}^d\Gamma_{i}$, where $\Gamma_i$ is the image of $\x_i(\O)$, for $i = 1,\dots,d$. Let $\rho(\bm{\xi})$ be the probability density function(p.d.f) of $\bm{\xi}$. Then the problem \eqref{sto_eq} can be restated as following: 
find $u:\bar{D}\times \Gamma\rightarrow \R$ such that $\rho$-a.e. for $\bm{\xi}\in \Gamma$, the equation
\begin{equation}\label{sto_eq1}
u_t(\mb{x},t;\bm{\xi}) = \mathcal{L}(\bm{x},t,\bm{\xi};u), \qquad \mb{x}\in D,
\end{equation}
with appropriate initial and boundary conditions. 

\subsection{Multi-element decomposition in random space}
To deal with irregularity (e.g. discontinuities) in the random space dependence we discretize the parametric space $\Gamma$ into a collection of non-overlapping hypercubes. Since any $d$-dimensional independent random vector can be transformed to a $d$-dimensional independent uniform random vector, without loss of generality, let $\bm{\xi}$ be a $d$-dimensional random vector defined on the random space $\Gamma$, where $\xi_i,i=1\cdots,d$ are independent identically distributed (i.i.d) uniform random variables defined on $[-1,1]$ with p.d.f $f_i=\frac{1}{2}$. Then $\Gamma = [-1,1]^d\subset\R^d$ with joint p.d.f. $\rho = \frac{1}{2^d}$. Following the notation in \cite{WanK_JCP05}, we briey introduce the decomposition
of the uniform random space (see \cite{WanK_JCP05} for more details). Let $\Gamma$ be decomposed into $N$ non-overlapping hypercubes as following:
\begin{eqnarray}\label{decom}
&&B_k = [a_1^k,b_1^k)\times[a_2^k,b_2^k)\times\cdots\times[a_d^k,b_d^k],\nonumber\\
&&\Gamma = \bigcup_{k=1}^N B_k,\\
&&B_i\cap B_j = \emptyset \quad \text{if } i\neq j,\nonumber
\end{eqnarray}
where $i,j,k = 1,2,\cdots,N$. Let $\chi_k,k=1,2,\cdots,N$ be the indicator random variables on each of the elements defined by
\[
\chi_k(\bm{\xi}) = \left\{\begin{array}{ll}1&\text{if }\bm{\xi}\in B_k,\\
0&\text{otherwise.}
\end{array}\right.
\]
For each random element $B^k$, the local random vector is defined by
\[
\bm{\zeta}^k = (\z_{1}^k,\z_{2}^k,\cdots,\z_d^{k})= \bm{\xi}|_{B^k},
\]
subject to the conditional p.d.f
\[
f_{\bm{\z}^k} = \frac{1}{2^d\Pb(\chi_k=1)}, \quad k=1,2,\cdots,N,
\]
where $\Pb(\chi_k=1) = \prod_{i=1}^d\frac{b_i-a_i}{2}$.
After that, we transfer each $\bm{\z}^k$ to a new random vector $\bm{\xi}^k$defined on $[-1,1]^d$ by a map $g_k$,
\[
g_k(\bm{\z}^k):\xi^k_i = \frac{b_i^k-a_i^k}{2}\z_i^k+\frac{b_i^k+a_i^k}{2},\quad i=1,2,\cdots,d.
\]
Thus, $\bm{\xi}^k$ is a random vector defined on $[-1,1]^d$ with constant p.d.f $f^k = \frac{1}{2^d}$.
With this decomposition, we can solve a system of differential equations with random input $\bm{\xi}$ by combining the local approximations via $\bm{\z}^k$ subject to a conditional p.d.f.
In practice if the solution $u(\bm{\xi})$ is locally approximated by $\hat{u}_k(\bm{\z}^k),k=1,2,\cdots,N$, then the $m$th moment of $u(\bm{\xi})$ on the entire random space can be obtained by
\begin{equation}\label{m_moment}
\begin{split}
\mu_m(u(\bm{\xi})) &= \int_{\Gamma}u^m(\bm{\xi})\frac{1}{2^d}d\bm{\xi}\\
&\approx\sum_{k=1}^N\Pb(\chi_k=1)\int_{B^k}\hat{u}_k^m(\bm{\z}^k)f_{\bm{\z}^k}d\bm{\z}^k\\
&=\sum_{k=1}^N\Pb(\chi_k=1)\E(\hat{u}^m(\mb{x},t;\cdot)|\bm{\xi}\in B^k) \\
&=\sum_{k=1}^N\Pb(\chi_k=1)\int_{\Gamma}\hat{u}_k^m(g_k^{-1}(\bm{\xi}^k))f^k(\bm{\xi}^k)d\bm{\xi}^k.
\end{split}
\end{equation}

\subsection{Localized stochastic discretization}
To solve the equation \eqref{sto_eq}, any stochastic method (such as Monte Carlo simulations(MCS), general polynomial chaos (gPC) Galerkin methods, stochastic collocation methods and etc.) can be applied in random space and result in a set of deterministic equations in the physical space which can be solved by standard numerical techniques. In the multi-element framework, MCS is straightforward, while both gPC Galerkin and stochastic collocation methods can be adapted to fit in. 

Next, we will briefly introduce the two most popular intrusive methods to deal with random inputs in one element.
\subsubsection{ME-gPC}
Without loss of generality, consider an orthonormal generalized polynomial chaos basis $\{\Phi_\i\}_{|\i|=0}^{\infty}$ spanning the space of second-order random processes on the probability space element $B^k$ ($\i = (i_1,\cdots,i_d)\in \mathbb{N}^d_0$ is a multi-index with $|\i|=i_1+\cdots+i_d$). The basis functions $\Phi^k_\i(\bm{\z}^k(\o))$ are polynomials of degree $|\i|$ with orthonormal relation
\begin{equation}\label{orth_basis}
\langle\Phi_\i^k,\Phi^k_{\j}\rangle_{B^k} = \delta_{\i\j},
\end{equation}
where $\delta_{\i\j}$ is the Kronecker delta and the inner product between two functions $g(\bm{\x})$ and $h(\bm{\x})$ is defined by
\begin{equation}\label{inner}
\langle g,h\rangle_{B^k}=\int_{B^k}g(\bm{\z}^k)h(\bm{\z}^k)f_{\bm{\z}^k}(\bm{\z}^k)d\bm{\z}^k.
\end{equation}

The solution $u$ is approximated by the truncated gPC expansion
\begin{equation}\label{gpc_solu}
u(\mathbf{x},t;\bm{\z}^k) = \sum_{|\i|=0}^p\hat{u}^k_\i(\mathbf{x},t)\Phi^k_\i(\bm{\z}^k(\o)).
\end{equation}
Substituting equation \eqref{gpc_solu} into the governing system \eqref{sto_eq}, we obtain the following system
\begin{equation}\label{sto_eq_gpc}
\sum_{|\i|=0}^p\frac{\partial\hat{u}^k_{\i}}{\partial t}\Phi^k_\i=\mathcal{L}(\mathbf{x},t,\bm{\z}^k;\sum_{|\i|=0}^p \hat{u}^k_\i\Phi^k_\i).
\end{equation}
By using the Galerkin projection of \eqref{sto_eq_gpc} onto each element with the orthonormal polynomial basis $\{\Phi^k_\i\}_{|\i|=0}^p$, we derive
\begin{equation}\label{sto_eq_galerkin}
\frac{\partial\hat{u}_{\i}^k}{\partial t}=\langle\mathcal{L}(\mathbf{x},t,\bm{\z}^k;\sum_{|\i|=0}^p \hat{u}^k_\i\Phi^k_\i),\Phi^k_\j\rangle_{B^k} ,\quad |\j| = 0,1,\cdots,p.
\end{equation}
This is a set of $n = \binom {d+p} d$ coupled deterministic equations of the random modes $\hat{u}^k_\i(\mathbf{x},t),|\i| = 0,1,\cdots,p.$ Techniques for deterministic equations can be implemented to solve this system of equations. Then the mean of the solution on element $B^k$ can be evaluated by
\begin{equation}\label{condi_m_moment_ga}
\begin{split}
\E(\hat{u}(\mb{x},t;\cdot)|\bm{\xi}\in B^k) &= \int_{B^k}\hat{u}(\mb{x},t;\bm{\xi}))f(\bm{\xi})d\bm{\xi}\\
&\approx \int_{B^k}\sum_{|\i|=0}^p\hat{u}^k_{\i}(\mb{x},t)\Phi_\i(\bm{\z}^k)f_{\bm{\z}^k}(\bm{\z}^k)d\bm{\z}^k.\\
\end{split}
\end{equation}
The $m$-th moment of the  approximated solution $\hat{u}(\mb{x},t;\bm{\x})$ over element $B^k$ can be evaluated similarly. The global $m$-th moment of  $\hat{u}(\mb{x},t;\bm{\x})$ can be obtained by \eqref{m_moment}.

\subsubsection{ME-PCM}
To avoid solving a large coupled system of deterministic equations one can use the multi-element polynomial collocation method (ME-PCM) \cite{Foo2008}. To implement ME-PCM, we need to specify the set of collocation points. Once the mesh discretization of $\Gamma$ is prescribed, a set of collocation points $\{\mb{z}_j^k\}_{j = 1,r}$ is identified in each element $B_k$, where $r$ is the number of points used. Usually, these points are chosen to be the points of a cubature rule on $B_k$ with integration weights $\{w_j^k\}_{j=1}^r$, $k = 1,\dots,N$. Here we consider full tensor products of Gauss quadrature points for $d=1,2$, and sparse grid points for $d\geq 3$.

At each collocation point $\mb{z}^k_j$, we find the solution of the deterministic problem
\begin{equation}\label{sto_eq_det}
u_t(\mb{x},t;\mb{z}_j^k) = \mathcal{L}(\mb{x},t,\mb{z}_{j}^k;u), \qquad \mb{x}\in D,
\end{equation}
with appropriate initial and boundary conditions.
Then for each element $B_k$, $k=1,\dots,N$, we can construct an approximate solution $\mathcal{I}^1_{B^k}u(\mb{x},t;\bm{\x})$ using the set of solutions on the collocation points. For instance the operator can be chosen to be the tensor product Lagrangian interpolant, i.e.,
\begin{equation}\label{eq_lag}
\mathcal{I}^1_{B^k}u(\mb{x},t;\bm{\x}) = \sum_{j=1}^{r} u(\mb{x},t;\mb{z}_{j}^k)\ell_{j}^k(\bm{\z}^k),
\end{equation}
where $\ell_{j}^k(\bm{\z}^k)$ is the Lagrange polynomial corresponding to the collocation points $\{\mb{z}_j^k\}_{j=1}^r$, refer \cite{Foo2008} for details. 
Another frequently used operator is the orthogonal polynomial interpolant (gPC interpolant). To do this interpolation, assume in element $B^k$, $u(\mb{x},t;\bm{\x})$ can be approximated by $\hat{u}(\mb{x},t;\bm{\x}) = \sum\limits_{|\i|=0}^p\hat{u}_{\i}^k(\mb{x},t)\Phi_{\i}^k(\bm{\z}^k)$, then by the cubature rule and the orthonomality of $\Phi_{\i}^k(\bm{\z}^k)$, we have 
\begin{equation}\label{def_ui_col}
\hat{u}_{\i}^k(t,\mb{x}) = \sum_{j=1}^r u(t,\mb{x};\mb{z}^k_j)\Phi^k_{\i}(\mb{z}^k_j)w_{j}^k, \quad |\i| = 0,\cdots,p.
\end{equation}
The gPC interpolation is constructed by
\begin{equation}\label{eq_gpc}
\mathcal{I}^2_{B^k} u(\mb{x},t;\bm{\x})  = \sum_{|\i|=0}^p\hat{u}_{\i}^k(\mb{x},t)\Phi_{\i}^k({\bm{\z}^k}).
\end{equation}
As a result, the global approximation is defined as:
\begin{equation}\label{eq_lag_glo}
\hat{u}(\mb{x},t;\bm{\x}) = \sum_{k=1}^N\mathcal{I}^i_{B^k}u(\mb{x},t;\bm{\x}) \chi_k(\bm{\x}), \quad \mb{x}\in D \quad \bm{\x}\in\Gamma.
\end{equation}
where $i = 1$ refers to Lagrangian interpolation and $i=2$ refers to gPC interpolation.

Considering the statistical properties of the solution $\hat{u}(\mb{x},t;\mb{z})$, once the m-th moment of $\hat{u}$ is calculated on each stochastic element the global $m$-th moment can be assembled by \eqref{m_moment}. The conditional $m$-th moment  on element $B^k$ can be evaluated via the cubature rule over $B^k$ by
\begin{equation}\label{condi_m_moment}
\begin{split}
\E(\hat{u}^m(\mb{x},t;\cdot)|\bm{\xi}(\o)\in B^k) &= \int_{B_k}\hat{u}_k^m(\mb{x},t;\bm{\z}^k)f_{\bm{\z}^k}d\bm{\z}^k\\
&\approx \sum_{j=1}^r u^m(\mb{x},t;\mb{z}_j^k)w_{j}^k.
\end{split}
\end{equation}
Again, we can obtain the global $m$-th moment by \eqref{m_moment}. For the properties of this discretization and error analysis please refer to \cite{Foo2008}.

\section{Adaptive mesh refinement}\label{sec:mesh_refinement}
The structure of the division of the random space as described in the previous section can be allowed to change with time. In particular, different regions of random space can require different resolutions. This situation requires an adaptive mesh refinement procedure. In \cite{Li2015} we introduced an adaptive  mesh refinement algorithm based on the reduced model of the full system. Once a reduced model is obtained, a quantity  guiding the refinement is calculated. Unfortunately, a good reduced model of a general system is seldom available in an explicit form. However, the idea to use the memory terms in reduced model inspires the quantity used in this paper for the refinement. 

To convey the main idea, we can focus on one random space element. Suppose that the truncated Galerkin projection of equation \eqref{sto_eq} with orthonormal basis is available, and assume $p$ is the highest degree of the polynomial basis in the truncated system, then we obtain the full model 
\begin{equation} \label{sto_eq_galerkin_p}
\frac{\partial\hat{u}_{\j}}{\partial t}=\langle\mathcal{L}(\mathbf{x},t,\bm{\x};\sum_{|\i|=0}^p \hat{u}_\i\Phi_\i),\Phi_\j\rangle ,\quad |\j| = 0,1,\cdots,p.
\end{equation}
For the reduced system, choose $p_0<p$ (for our numerical examples we have chosen the ratio $p_0/p$ to be around $1/2$). 


Then the truncated system up to $p_0$ can be written as:
\begin{equation}\label{sto_eq_galerkin_p0}
\frac{\partial\t{u}_{\j}}{\partial t}=\langle\mathcal{L}(\mathbf{x},t,\bm{\x};\sum_{|\i|=0}^{p_0} \t{u}_\i\Phi_\i),\Phi_\j\rangle ,\quad |\j| = 0,1,\cdots,p_0.
\end{equation}
Let $$\mathbf{E}_{p_0}(\hat{u}) := \sum\limits_{|\j|=0}^{p_0}\hat{u} ^2_{\j} = \|\sum\limits_{|\j|=0}^{p_0} \hat{u}_{\j}\Phi_{\j}(\bm{\x})\|^2_{\Gamma},$$ and $$\mathbf{E}_{p_0}(\t{u}) := \sum\limits_{|\j|=0}^{p_0}\t{u} ^2_{\j} = \|\sum\limits_{|\j|=0}^{p_0}\t{u}_{\j}\Phi_{\j}(\bm{\x})\|_{\Gamma}^2.$$ $\mb{E}_{n}(u)$ can be interpreted as the $L^2$ norm of the projection of $u$ onto the space spanned by $\{\Phi_{\j}(\bm{\xi})\}_{|\j|=0}^n$. Then the rate of change of $\mathbf{E}_{p_0}(\hat{u})-\mathbf{E}_{p_0}(\t{u})$ denoted by $Q$ is obtained as:
\begin{equation}\label{rate_of_change1}
\begin{split}
Q:=&\frac{d(\mb{E}_{p_0}(\hat{u})-\mb{E}_{p_0}(\t{u}))}{dt} =2\sum_{|\j|=0}^{p_0}\frac{\partial \hat{u}_{\j}}{\partial t} \hat{u}_{\j}-2\sum_{|\j|=0}^{p_0}\frac{\partial \t{u}_{\j}}{\partial t} \t{u}_{\j}\\
=&2 \sum_{|\j|=0}^{p_0} \langle\mathcal{L}(\mathbf{x},t,\bm{\x};\sum_{|\i|=0}^{p} \hat{u}_\i\Phi_\i),\Phi_\j\rangle {\hat{u}}_{\j}-2\sum_{|\j|=0}^{p_0} \langle\mathcal{L}(\mathbf{x},t,\bm{\x};\sum_{|\i|=0}^{p_0} \t{u}_\i\Phi_\i),\Phi_\j\rangle{\t{u}}_{\j}.
\end{split}
\end{equation}
$Q$ measures the energy transfer between the full system (truncated up to order $p$) and the reduced system (truncated up to order $p_0$). Let $\mb{Q} = \int_{D}|Q|d\mb{x}$.
When the quantity $\mb{Q}$ multiplied by the volume of a specific element exceeds a user-prescribed tolerance then the element is split in two parts i.e., mesh refinement takes place (see also \cite{Li2015}).

If the random space has $d$ dimensions with $d \geq 2$ then we need to decide not only when to refine but also in which directions. This can be decided in different ways and here we offer two possibilities. 

The first is to consider the contribution of the highest degree ($p_0$) basis in that dimension in the reduced model. In order to implement we define
\begin{equation}\label{s1_i}
\begin{split}
s^1_i =&\int_D \big|2\langle \mathcal{L}(\mb{x},t,\bm{\x};\sum_{|\i|=0}^p\hat{u}_{\i}\Phi_{\i}),\Phi_{p_0\mb{e}_i}\rangle \hat{u}_{p_0\mb{e}_i}\\
&-2\langle \mathcal{L}(\mb{x},t,\bm{\x};\sum_{|\i|=0}^{p_0}\t{u}_{\i}\Phi_{\i}),\Phi_{p_0\mb{e}_i}\rangle \t{u}_{p_0\mb{e}_i}\big|d\mb{x},
\end{split}
\end{equation}
where $i$ refers to the $i$-th dimension, and $\mb{e}_i$ is the index vector with $i$-th dimension 1 others 0. Compared with \eqref{rate_of_change1}, $s^1_i$ is the contribution in degree $p_0$ basis in the $i$-th dimension. 

The second criterion for directional refinement is related to the total variance in the $i$-th direction. We define
\begin{equation}\label{s2_i}
\begin{split}
s^2_i =& \int_{D}\big|2\sum_{n=1}^{p_0}\langle \mathcal{L}(\mb{x},t,\mb{\x};\sum_{|\i|= 0}^{p}\hat{u}_{\i}\Phi_{\i}),\Phi_{n\mb{e}_i}\rangle\hat{u}_{n\mb{e}_i}\\\
&-2\sum_{n=1}^{p_0}\langle \mathcal{L}(\mb{x},t,\mb{\x};\sum_{|\i|= 0}^{p}\t{u}_{\i}\Phi_{\i}),\Phi_{n\mb{e}_i}\rangle\t{u}_{n\mb{e}_i}\big|d\mb{x}.\\
\end{split}
\end{equation}
The quantity $s^2_i$ is the contribution to $Q$ from all the degrees $1$ to $p_0$ in the $i$-th dimension. 

We now have all the ingredients needed to present our mesh refinement algorithm:

\vspace{0.5cm}

\begin{enumerate}
\item[]{\bf Mesh refinement algorithm}
\item[Step 1] Choose values $TOL_1>0$ and $TOL_2>0$ for the tolerances.
\item[Step 2] Mesh refinement:\\
For time step $t \leftarrow 1, \cdots, N_T$ \\
\begin{itemize}
\item[] Loop over all elements:
\begin{itemize}
\item[] On the $k$-th element $B_k$, update the modes for $B_k$\\
 If {$\mb{Q}\Pr(B_k)\geq TOL_1$,} and if $d=1$, split element $B^k$ and update system, otherwise loop over all dimensions:
\begin{itemize}
\item[] If {$s^1_i ( \textrm{ or } s^2_i, \textrm{respectively})\geq TOL_2\cdot \max_{j=1,\cdots, d} s^1_j ( \textrm{ or } s^2_j, \textrm{respectively})$,}
\begin{itemize}
\item[] split the element $B_k$ in two equal parts along the $i$th dimension and update the system
\end{itemize}
\item[] End if
\end{itemize}
\item[] End if
\item[] Update the information of the new elements
\end{itemize}
\item[] End loop
\end{itemize}
\end{enumerate}
We have to make several remarks about the implementation details of the algorithm. 

First, in \eqref{rate_of_change1} we obtain an expression for the rate of change($Q$) which does not involve temporal derivatives. There are two different ways to obtain $Q$ which are based on the method used to solve \eqref{sto_eq}. If the gPC Galerkin method is employed to solve \eqref{sto_eq}, then at every time step the coefficients in the gPC expansion are available. In this case, $Q$ is computed directly via \eqref{rate_of_change1}. If the gPC collocation method is utilized to solve \eqref{sto_eq}, then extra steps to evaluate the coefficients of the gPC expansion (see \eqref{def_ui_col}) need to taken. However, since when applying collocation method the resulting system is decoupled, the collocation approach is usually more efficient than the Galerkin approach. 

Secondly, whether one uses a Galerkin or a collocation method, there are two ways to compute $Q$ then $\mb{Q}.$ If we are employing a Galerkin method, then the first way is to evolve, for each element $B_k$, both the full and reduced systems of equations i.e., solve equations \eqref{sto_eq_galerkin_p} and \eqref{sto_eq_galerkin_p0}. On the other hand, If we are employing a collocation method, then for each element $B_k$ we use the tensor mesh with $r$ points in each dimension and evolve the system on this mesh. Afterwards, we use a cubature rule to get the coefficients $\hat{u}_{\i}^k$, $|\i| = 0,\dots,p$ via \eqref{def_ui_col}. Here $r$ depends on the highest order of polynomial representation of each random variable in the RHS of \eqref{sto_eq_galerkin}. To get the coefficients of the reduced system in element $B_k$, we use another tensor mesh with $r_0$ points in each dimension and evolve the system, then evaluate the coefficient $\t{u}^k_{\i}$ for polynomial bases up to order $p_0$. For both approaches, we obtain $Q$ by \eqref{rate_of_change1}.  If the explicit expressions for both full system and reduced systems (depending only on the coefficients of the expansion in the orthonormal basis) are available, then only $p+1$ collocation points in each dimension for full system and $p_0+1$ collocation points for the reduced system are required. If we are employing a Galerkin method, the second way to obtain the quantity $Q$ is to evolve {\it only} the full system \eqref{sto_eq_galerkin}, and let $\t{u}_{\i}^k = \hat{u}_{\i}^k$, for $ |\i| = 0,\dots,p_0$. If we are using a collocation method then we only evolve the system on a tensor mesh with $r$ points in each dimension. Evaluate the coefficients for $\hat{u}_{\i}^k$, $|\i|=0,\dots,p$, then follow the steps as in the Galerkin case to get $Q$. 

The third remark concerns the generation of values in the newly created elements after splitting of the original element. For the Galerkin approach, when the element $B^k$ is split into two elements $B^{k_1}$ and $B^{k_2}$, two new random variables $\z^{k_1}$ and $\z^{k_2}$ are generated according to the procedure in \cite{WanK_JCP05}. For the collocation approach, after the refinement a new mesh is generated. Let $\{\mb{z}^{k_1}_j\}_{j=1}^{r}$ and $\{\mb{z}^{k_2}_j\}_{j=1}^r$ be the collocation points in the newly generated elements $B^{k_1}$ and $B^{k_2}$, respectively. We can use interpolation (see \eqref{eq_lag} or \eqref{eq_gpc}) to obtain the system values $\hat{u}(\mb{x},\tau;\mb{z}^{k_1}_j)$ and $\hat{u}(\mb{x},\tau;\mb{z}^{k_2}_j)$, $j = 1,\dots,r$  (here $\tau$ is the time when the refinement occurs). For the numerical examples below, we use gPC interpolation to update the system values on the new collocation points, since the gPC expansion is already available after evaluating $\mb{Q}$. 

Our final remark is very important because the performance of the mesh refinement procedure relies heavily on it. For many cases, starting with only one element in random space is not enough, even if one uses high order gPC. If one insists on starting with one element, then the new points generated after the mesh refinement are not well placed to guarantee the accurate evolution of the solution. This is to be expected since an under-resolved representation of the random space initially will lead to erroneous identification of the regions that need to be refined. To avoid this one must have an adequate number of elements initially. If there is no other way to estimate how many elements are needed initially, one can decide by running several realizations with a different initial number of elements and monitoring the convergence of the results as the initial number of elements is increased.


\section{Numerical Examples} \label{sec:examples} We have tested our adaptive mesh refinement algorithm both with a stochastic Galerkin method (AMR-GA) and with a stochastic collocation method (AMR-CO) on variety of problems of increasing complexity. In particular, we have applied it to a single linear ODE with a one-dimensional random parameter, the nonlinear Kraichnan-Orszag three-mode system with random initial conditions, as well as the 1D Kuramoto-Sivashinsky (K-S) equation with random bifurcation parameter.

\subsection{One-dimensional linear ODE}
We begin by considering the simple ODE \cite{WanK_JCP05}
\begin{equation}\label{ex:ODE}
\frac{du}{dt} = -\kappa(\omega)u,\quad u(0;\omega) = u_0,
\end{equation}
where $\kappa(\omega)\sim U(-1,1)$. The exact solution of this equation is
\begin{equation}\label{ex:ODE_sol}
u(t,\omega) = u_0e^{-\kappa(\omega)t}.
\end{equation}
The statistical mean of the solution is $$\mu(u(t;\o)) = \{\begin{array}{ll} \frac{u_0}{2t}(e^{t}-e^{-t}),& t>0,\\u_0,&t=0,\end{array}$$ and the variance is $$\sigma^2(u(t;\o))=\{\begin{array}{ll}\frac{u_0^2}{4t}(e^{2t}-e^{-2t})-\frac{u_0^2}{4t^2}(e^{2t}+e^{-2t}-2),&t>0,\\0,& t=0.\end{array}$$
We will focus on the implementation of the adaptive mesh refinement with gPC collocation approach without an explicit expression of Galerkin projection. Also, we will evolve only the full system and use that to compute all the necessary quantities to decide when and where to refine. Since the random input is uniformly distributed, we choose Gaussian quadrature points as the collocation points. Let $p$ be the order of the full system and  $p_0 = \lceil \frac{p+1}{2}\rceil$ the order of the reduced system. From Gaussian quadrature we know that if the gPC expansion of the solution is up to order $p,$ then the highest polynomial order in random variable $z$ of RHS in \eqref{ODE_cola} for which the quadrature is exact equals $2p+1.$ This means that in order to evaluate exactly the integral of these polynomials we need $p+1$ collocation points in each dimension of each element. Let $\{z^{k}_i\}_{i=1}^{p+1}$ be the collocation points in element $B^k$, and $\{w^k_i\}_{i=1}^{p+1}$ are the Gaussian quadrature weights respectively. Also, let ${u}(\tau,z^k_i)$, $i = 1,\dots,p+1$ be the solutions to \eqref{ex:ODE} at points $\{z^{k}_i\}_{i=1}^{p+1}$ at time $\tau.$ Let 
\begin{equation}\label{sol_gar}
\hat{u}^k(t,z)= \sum\limits_{i=0}^p\hat{u}^k_i(t)\Phi_i(z),
\end{equation}
where $\{\Phi_i\}$ are the orthonormal Legendre 
polynomial basis functions. The coefficients $\hat{u}_i^k(\tau)$, $i = 1,\dots,p.$ can be obtained through \eqref{def_ui_col}. To estimate $\mb{Q}$, first we substitute \eqref{sol_gar} to left hand side(LHS) of \eqref{ex:ODE} and we obtain:
\begin{equation}\label{ODE_gar}
\sum_{i=0}^p\frac{d \hat{u}_i^k}{d t} \Phi_i(z)= - z u(t,z), \quad i = 0,\cdots,p.
\end{equation}
By the orthonomality of $\Phi_i$, $i=0,\dots,p$, we have the following expression:
\begin{equation}\label{ODE_cola}
\frac{d \hat{u}_i^k}{d t} = -\int_{B^k}zu(t,z)\Phi_i(z)dz, \quad i = 0,\dots,p.
\end{equation} 
By the quadrature rule, at time $\tau$, \eqref{ODE_cola} becomes
\begin{equation}\label{ODE_colb}
\frac{d \hat{u}_i^k}{d t} \Big|_{t=\tau}= -\sum_{j=1}^{p+1}z^k_j u(\tau,z^k_j)\Phi_i(z_j^k)w^k_j, \quad i = 0,\dots,p.
\end{equation}
Immediately, we have:
\begin{equation}\label{Q_a}
\frac{d}{dt}\sum_{i=0}^{p_0} (\hat{u}^k_i)^2\Big|_{t=\tau} = 2\sum_{i=0}^{p_0}\hat{u}^k_i\frac{d \hat{u}^k_i}{dt} = -2\sum_{i=0}^{p_0}\sum_{j=1}^{p+1}\hat{u}^k_iz_j^ku(\tau,z_j^k)\Phi_i(z_j^k)w_j^k.
\end{equation}
As mentioned in the previous paragraph we do not evolve the reduced truncated system up to order $p_0.$ Instead, we let $\t{u}_i^k(\tau) = \hat{u}_i^k(\tau)$, for $i=0,\dots,p_0$. To evaluate the RHS of \eqref{sto_eq_galerkin_p0}, we substitute the approximation $u(\tau,z_j^k) = \sum\limits_{i=0}^{p_0}\t{u}_{i}^k\Phi_{i}(z_j^k)$, $j=1,\dots,p+1$ and apply the quadrature rule for the integral. We find
\begin{equation}\label{Q_b}
\begin{split}
\frac{d}{dt}\sum_{i=0}^{p_0} (\t{u}^k_i)^2\Big|_{t=\tau} &= 2\sum_{i=0}^{p_0}\t{u}^k_i\frac{d \t{u}^k_i}{dt} \\
&= -2\sum_{i=0}^{p_0}\sum_{j=1}^{p+1}\t{u}^k_iz_j^k\sum_{s=0}^{p_0}\t{u}_s^k(\tau)\Phi_s(z_j^k)\Phi_i(z_j^k)w_j^k.
\end{split}
\end{equation}
By subtracting \eqref{Q_b} from \eqref{Q_a} we obtain $Q(\tau)$ and $\mb{Q} (\tau)= |Q(\tau)|.$ If $\mb{Q}(\tau)$ is greater than $TOL1$, we perform the refinement and use interpolation to estimate the values at the new collocation points. 

\begin{figure}[htbp]
   \centerline{
   \psfig{file=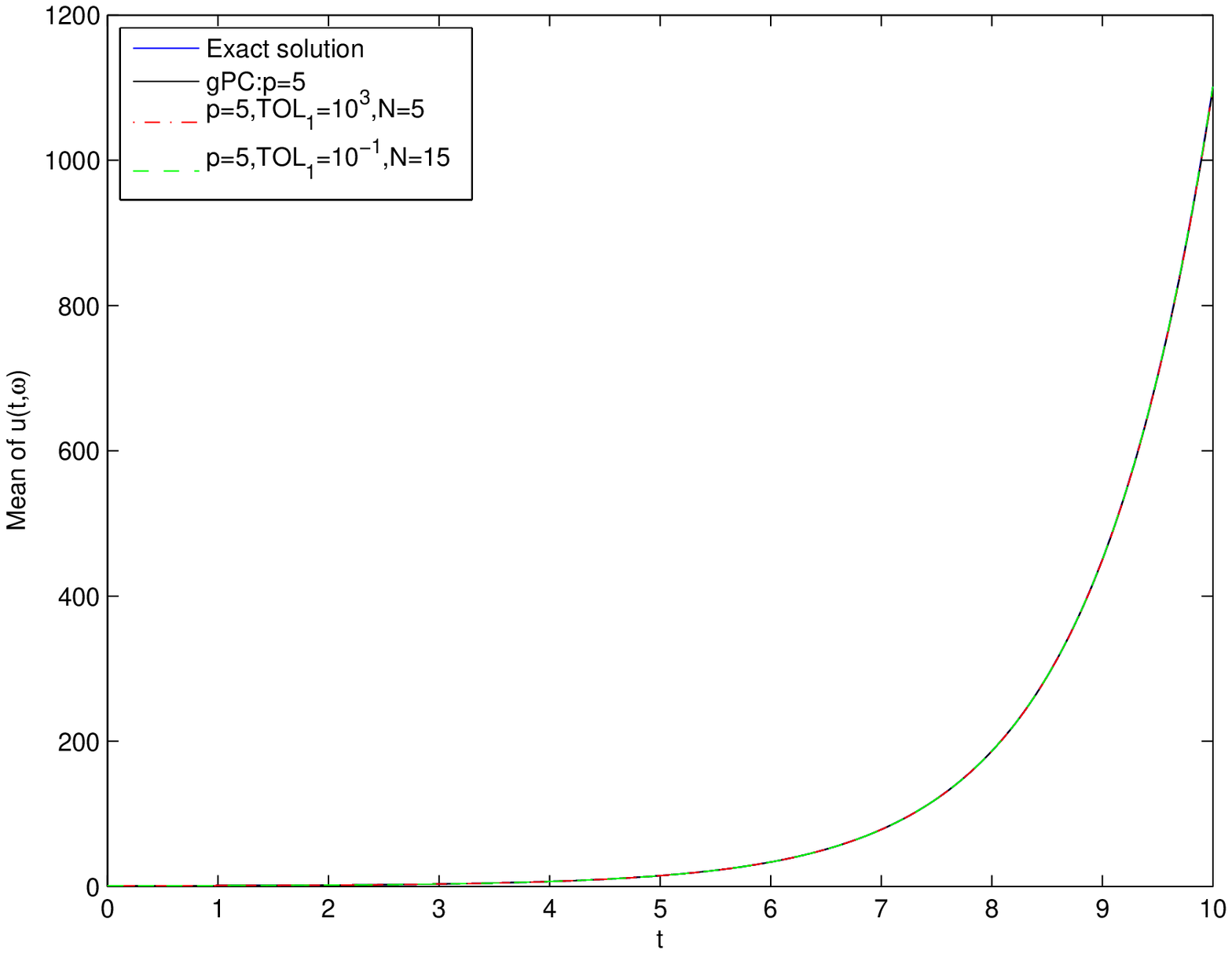,width=7cm}
   \psfig{file=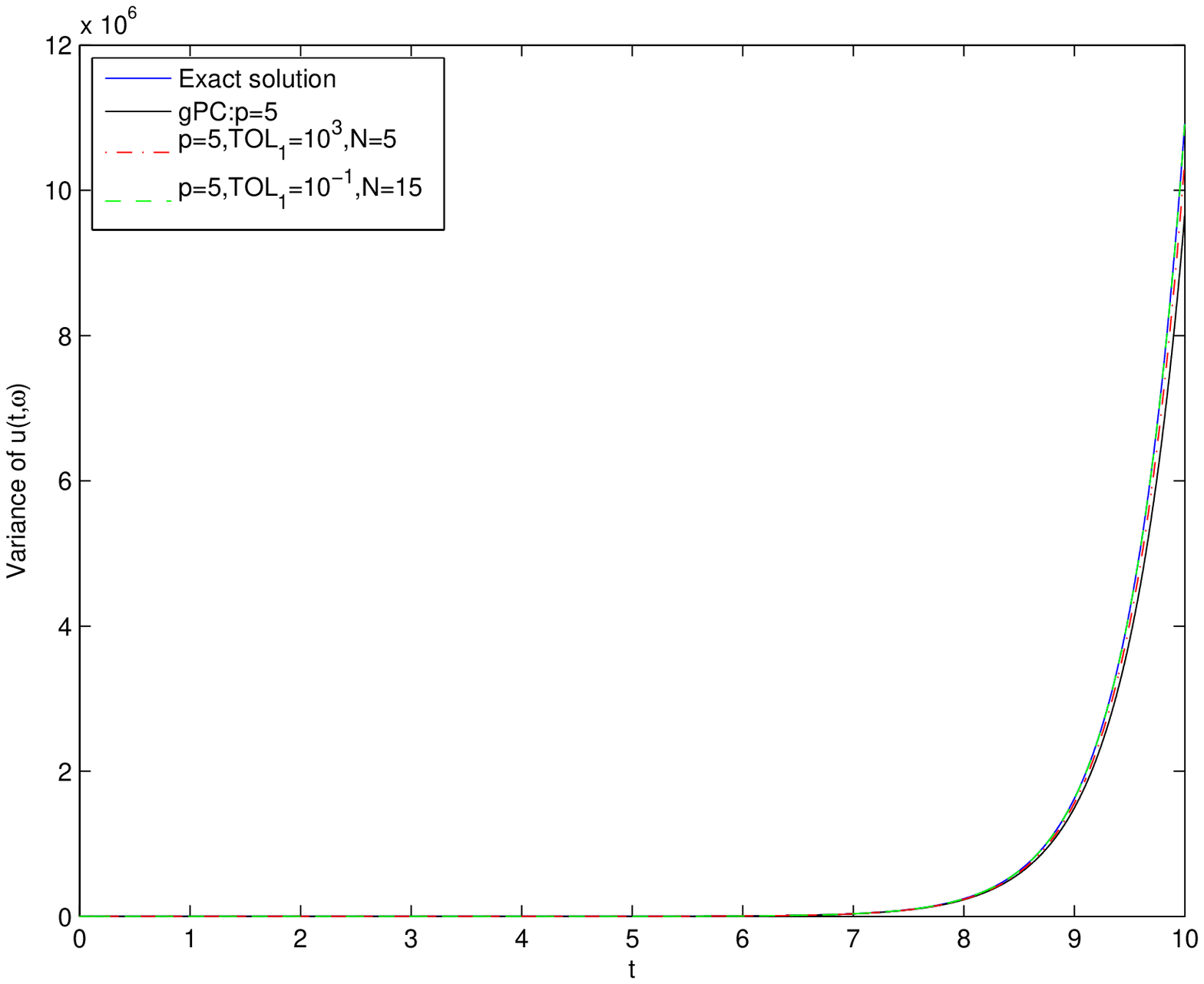,width=7cm}
  }
\caption{Evolution of mean of $u(t;\o)$(left) and evolution of variance of $u(t;\o)$(right) for the simple ODE.
  }
\label{fig:ode}
\end{figure}

\begin{figure}[htbp]
   \centerline{
   \psfig{file=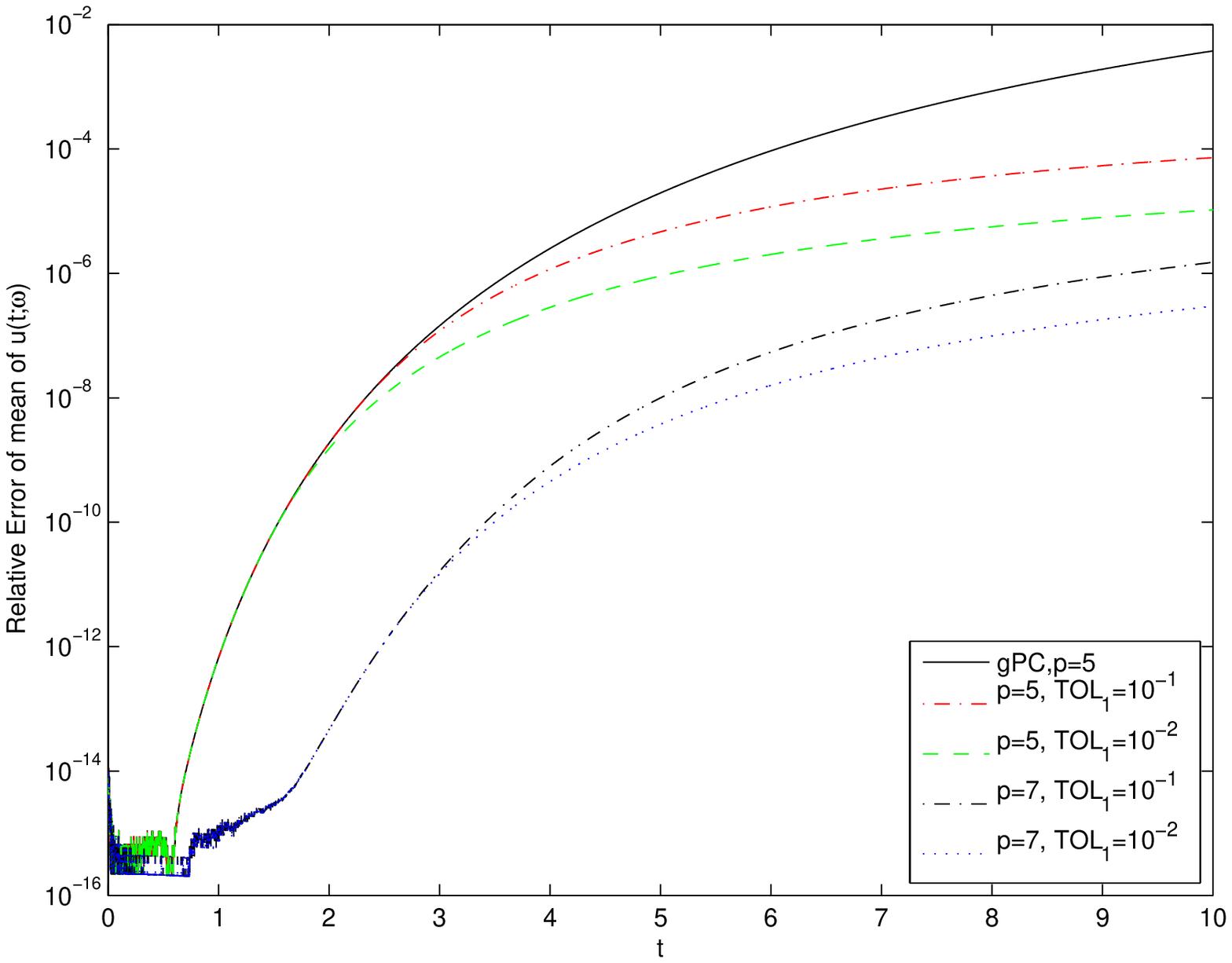,width=7cm}
   \psfig{file=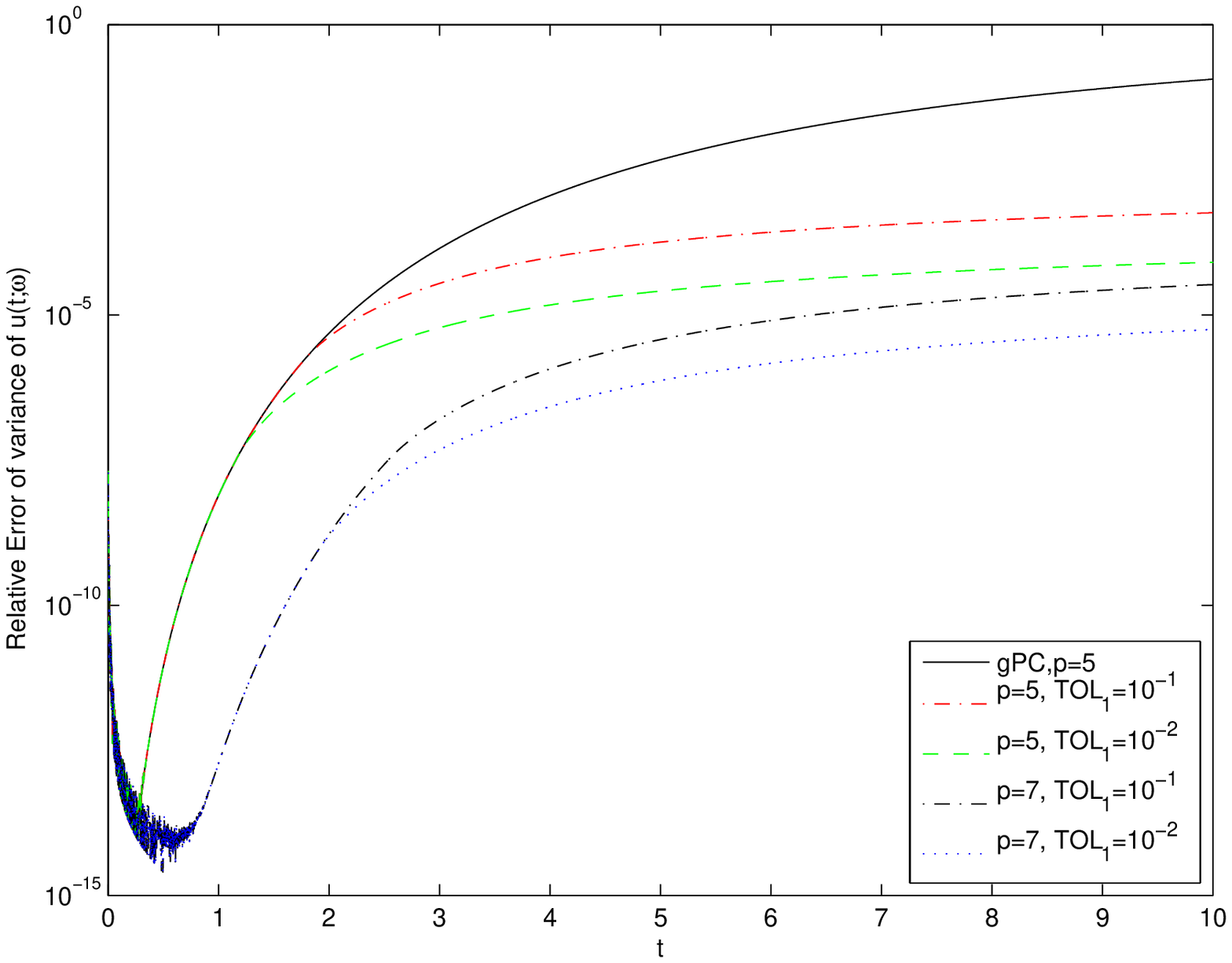,width=7cm}
  }
\caption{Evolution of relative error of mean of $u(t;\o)$ (left) and evolution of relative error of variance of $u(t;\o)$ (right) for the simple ODE.
  }
\label{fig:ode_er}
\end{figure}

\begin{figure}[htbp]
   \centerline{
   \psfig{file=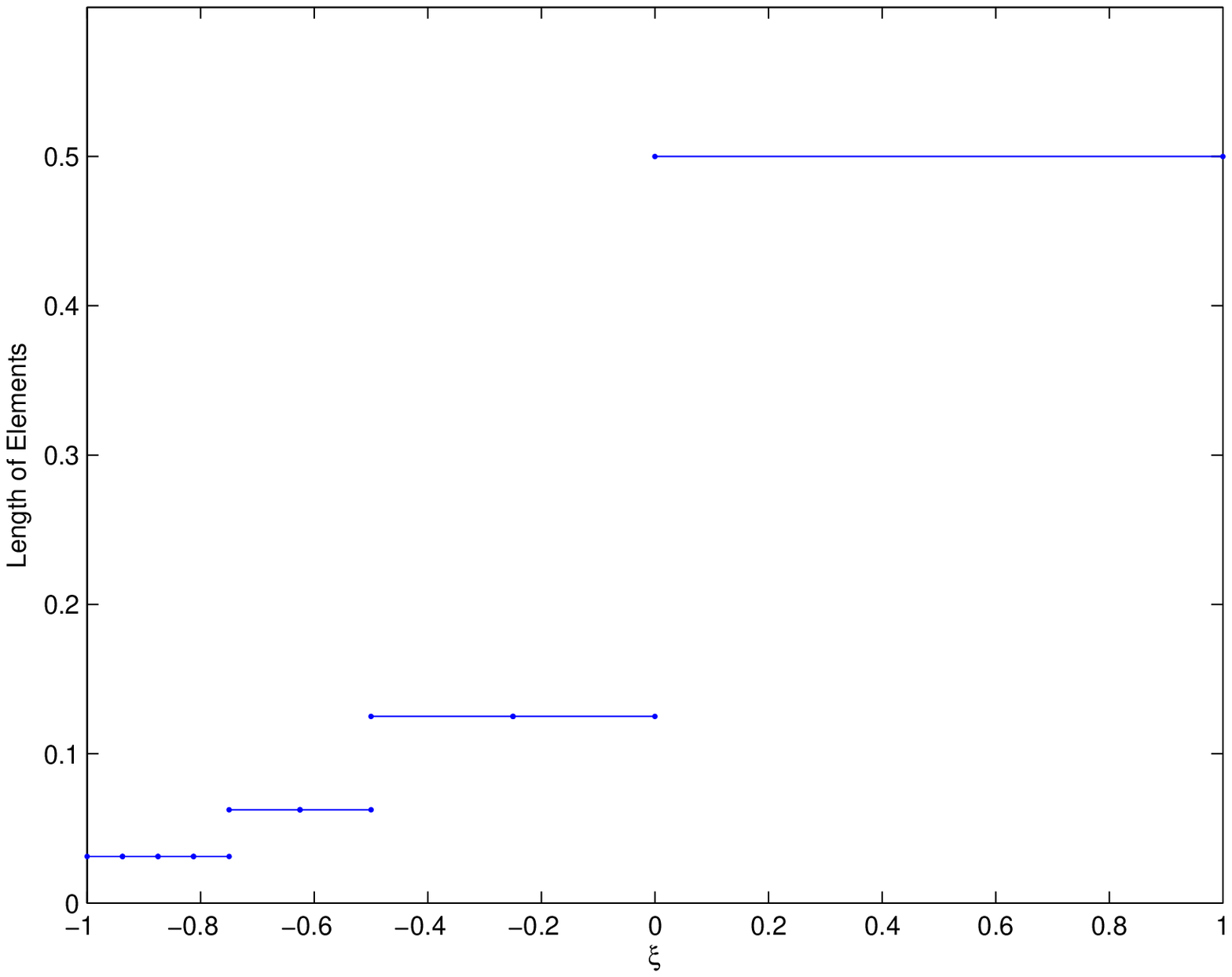,width=7cm}
   \psfig{file=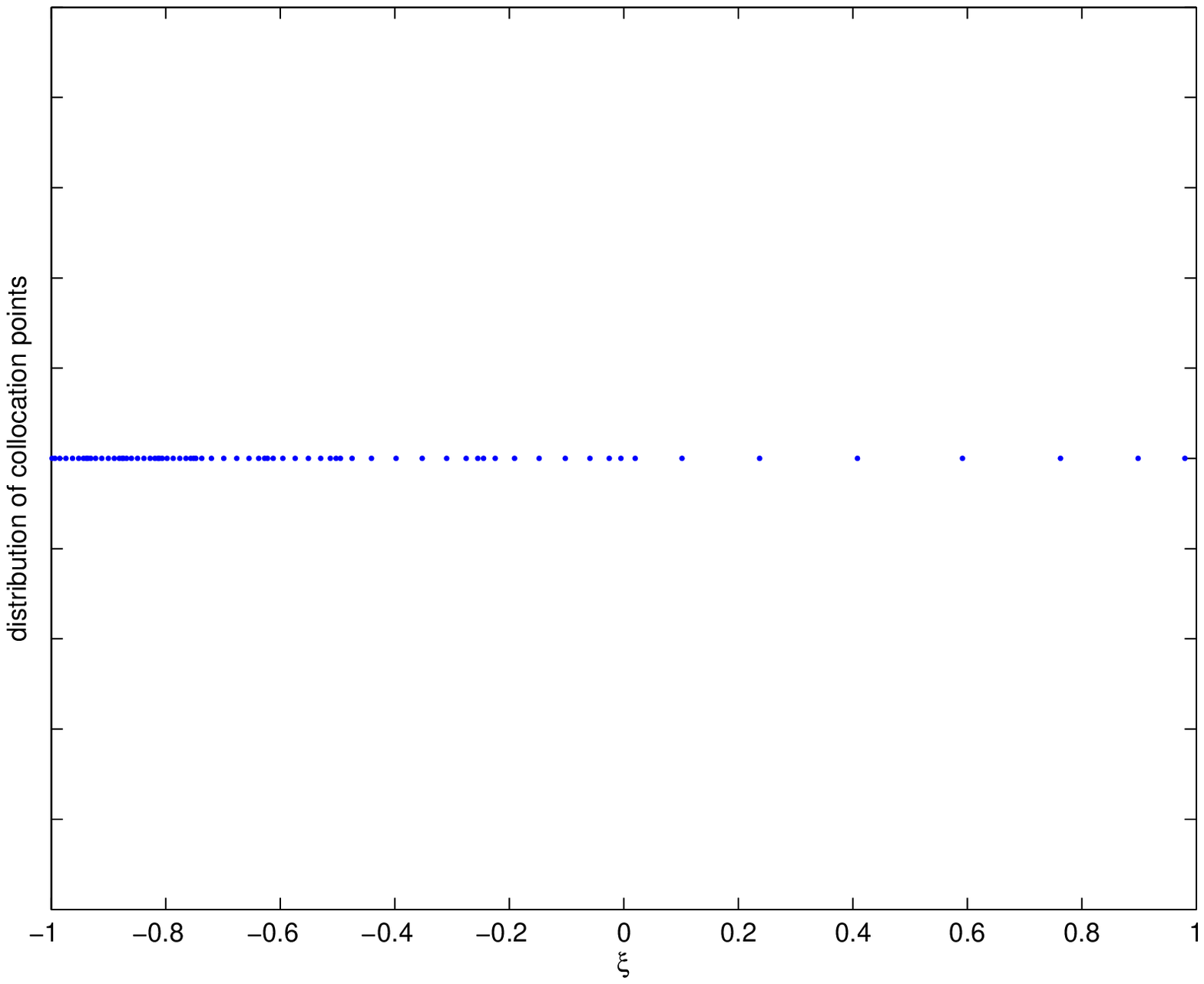,width=7cm}
  }
\caption{Result of adaptive mesh refinement for $p=7$, $TOL_1=10^{-1}$: length of the elements at $t = 10$ (left) and distribution of the collocation points at $t=10$ (right).
  }
\label{fig:ode_mesh}
\end{figure}

\begin{figure}[htbp]
   \centerline{
   \psfig{file=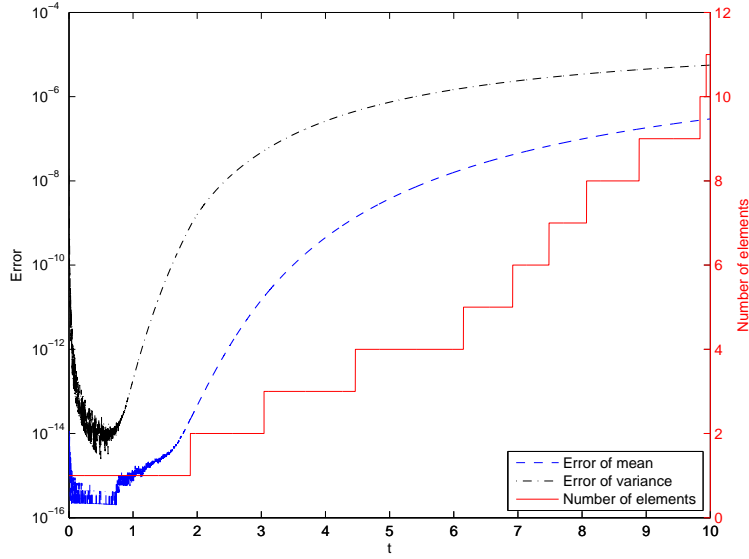,width=11cm}
  }
\caption{Evolutions of mean error and variance error compared with adaptive meshes for the simple ODE when  $p=7$, $TOL_1=10^{-2}$.
  }
\label{fig:ode_error_mesh}
\end{figure}

We study the evolution of the mean and the variance of the solution to \eqref{ex:ODE}. There is no discontinuity involved in this problem so we can start the evolution with only one element. However as time increases the lower ordered gPC solution (via Galerkin or collocation approach) starts to deviate from the exact solution. One way to keep the error under control is to increase the order of the gPC expansion in Galerkin approach or increase the number of collocation points in the collocation approach. For the adaptive mesh refinement algorithm we have chosen a low order gPC expansion for each element (fewer number of collocation points respectively). For the particular example, the result from the Galerkin method is the same as that from collocation method. Here we present only results for the collocation approach. 
Figure \ref{fig:ode} shows the curves of the mean (left) and variance (right) via various methods. The results from the refinement algorithm almost reproduce the exact solution while the global gPC solution starts departing from the exact solution as time evolves. We choose the mean and variance of the exact solution as references to study the relative error of each algorithm. We define
$$\textrm{Error of mean} = |\frac{\mu(u(t;\o))-\mu(\hat{u}(t;\o))}{\mu(u(t;\o))}|,$$
$$\quad\textrm{Error of variance} = |\frac{\sigma^2(u(t;\o))-\sigma^2(\hat{u}(t;\o))}{\sigma^2(u(t;\o))}|.$$
In Figure \ref{fig:ode_er} we present the evolutions of the error for the gPC solution with order 5, for the adaptive Galerkin approach of order 5 and the adaptive collocation approach of order 5 for different values of the accuracy control parameters. The maximum relative errors for the mean and the variance resulting from the adaptive mesh refinement with different orders are presented in Table \ref{tab:ode}. For reasons of comparison we include the relative error from the gPC solution of order 5. As expected the higher ordered models require fewer elements while at the same time they maintain higher accuracy. The adaptive meshes at $t=10$ with different accuracy control values $TOL_1=10^{-1}$ and $TOL_1 = 10^{-2}$ are demonstrated in Figure \ref{fig:ode_mesh}. The elements on the left end are smaller than those on the right end. This is consistent with the fact that the rate of change of $u(t;\o)$ is larger on the left end of the random space (this is because for negative values of the parameter the solution blows up exponentially and thus is more sensitive to uncertainty). The evolution of the error, of the number of elements and of the error for the variance are shown in Figure \ref{fig:ode_error_mesh}.

\begin{table} [htbp]
\begin{center}
\begin{tabular}{llccc}
\toprule
        & N & Error of $\mu(u)$ & Error of $var(u)$\\
\midrule
gPC, $p=5$& $1$ & $3.8e-3$ & $1.1e-1$\\
\midrule
$TOL_1=10^{-1}$\\
$p=5$ & $15$ & $7.3e-5$ & $5.7e-4$\\
$p=7$ & $9$ & $1.5e-6$ & $3.3e-5$\\
\midrule
$TOL_1=10^{-2}$\\
$p=5$ & $19$ & $1.0e-5$ & $8.0e-5$ \\
$p=7$ & $11$ & $3.0e-7$ & $5.6e-6$ \\
\bottomrule
\end{tabular}
\caption{Maximum of relative error of mean and variance of solution to the simple ODE when $t\in(0,10]$}
\label{tab:ode}
\end{center}
\end{table}

\subsection{The Kraichnan-Orszag three-mode system}
Consider the following system obtained by a linear transformation performed on the original Kraichnan-Orszag (K-O) three-mode system \cite{WanK_JCP05,Foo2008,Ma2009},
\begin{equation}\label{ex:KO}
\begin{split}
\frac{dy_1}{dt} &= y_1y_3,\\
\frac{dy_2}{dt} &= -y_2y_3,\\
\frac{dy_3}{dt} &= -y_1^2+y_3^2,
\end{split}
\end{equation}
subject to initial conditions,
\begin{equation}\label{ex:KO_ini}
y_1(0) = y_1(0;\omega),\quad y_2(0) = y_2(0;\omega),\quad y_3(0) = y_3(0;\omega).
\end{equation}
The solution of this system changes qualitatively depending on the value of the initial conditions $y_1(0)$ and $y_2(0)$. For certain initial conditions the solution of the problem is periodic, and the period goes to infinity
as $y_1(0)$ and $y_2(0)$  approach to the planes $y_1 = 0$ and $y_2 = 0$. Thus a discontinuity in the initial condition dependence exists when the initial conditions are allowed to be on both sides of these two planes. It was shown that a global Wiener-Hermite expansion does not faithfully represent the dynamics of the system when the random inputs are Gaussian variables in \cite{Orszag1967}. Mesh refinement algorithms  \cite{WanK_JCP05,Foo2008} and the adaptive hierarchical sparse grid collocation \cite{Ma2009} can efficiently quantify the uncertainty of the system when the random inputs are uniform random variables and are in a range which includes the discontinuity value. We consider the case with one random input, two random inputs and three random inputs respectively and the random initial conditions coming from a uniform distribution. We present results for both the Galerkin and collocation approaches. The associated system of ordinary differential equations is solved using the standard fourth-order Runge-Kutta scheme with a timestep $\Delta t = 0.01$.

\subsubsection{One-dimensional random input}
We choose the initial conditions
\begin{equation}\label{ex:KO_ini_1d}
y_1(0;\o) = 1,\quad y_2(0;\o) = 0.1\xi(\o),\quad y_3(0;\omega)=0,
\end{equation}
where $\xi\sim U[-1,1]$. In this case the discontinuity point $y_2=0$ is in the random input space.

\begin{figure}[htbp]
   \centerline{
   \psfig{file=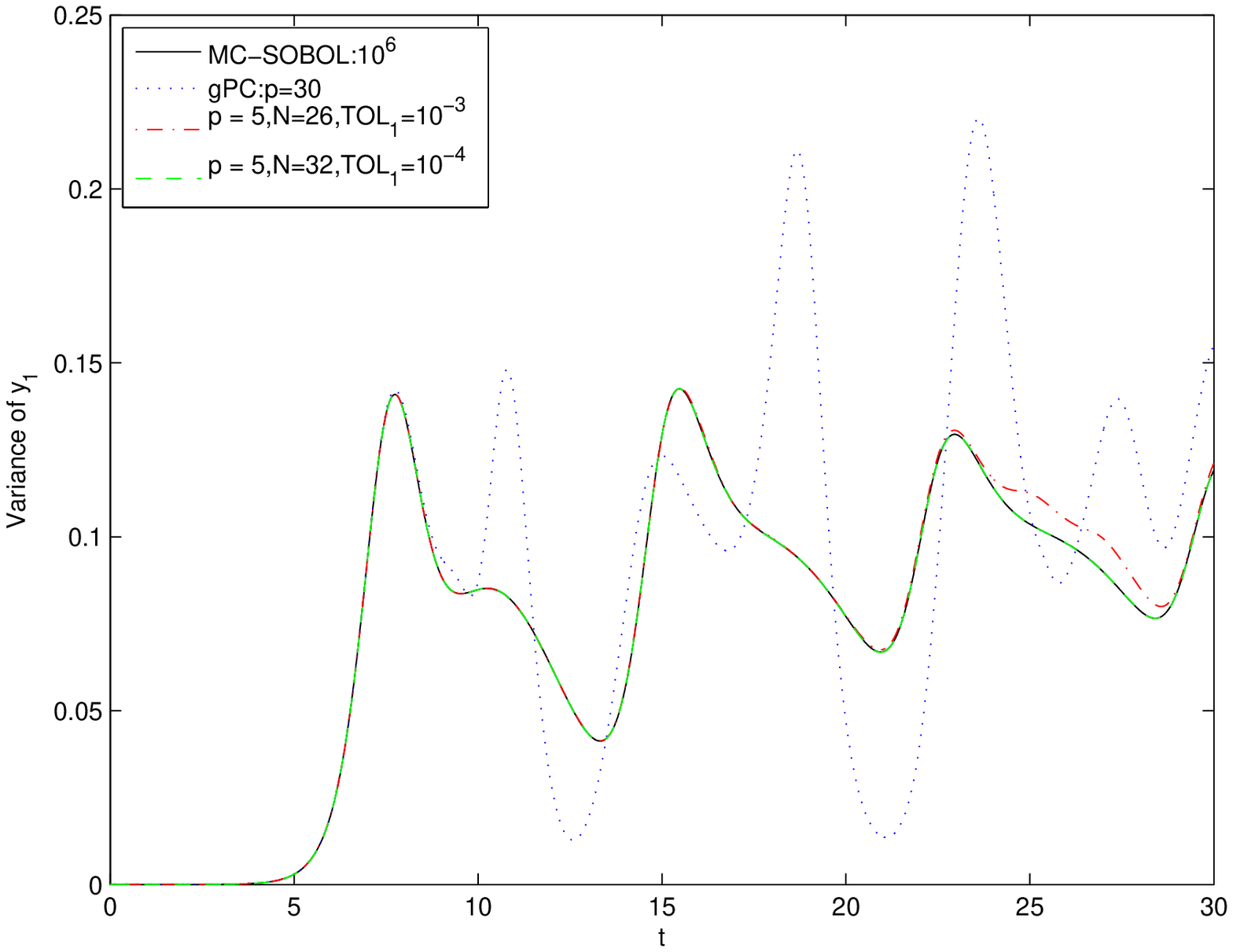,width=7cm}
   \psfig{file=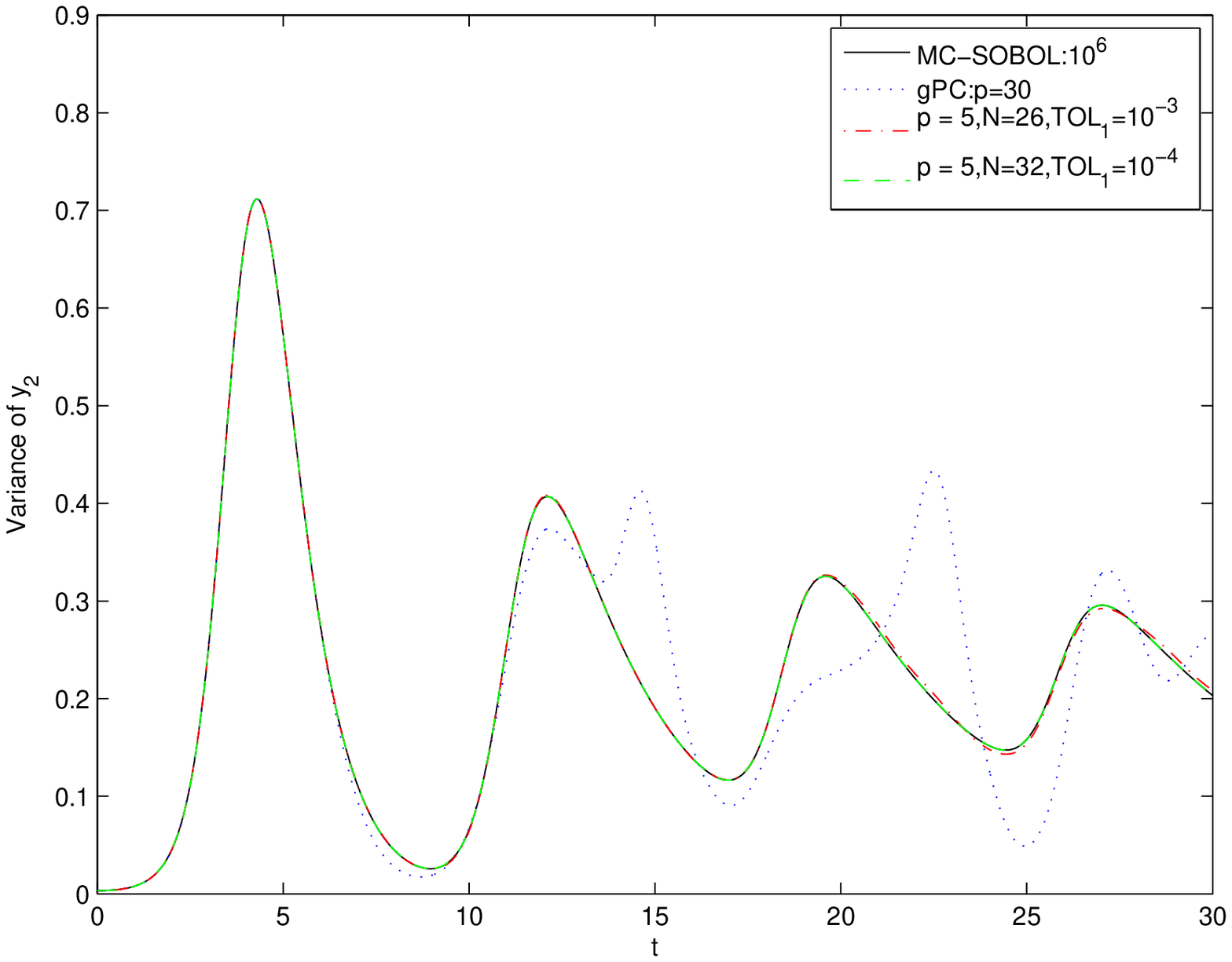,width=7cm}
  }
\caption{Evolution of the variance of $y_1$(left) and evolution of the variance of $y_2$(right) for the Kraichnan-Orszag three-mode system with 1D random initial condition.
  }
\label{fig:KO_y1y2_1d}
\end{figure}

We study the evolution of each variable $y_i$, $i=1,2,3$ on the time interval $[0,30]$ (short time behavior). Figure \ref{fig:KO_y1y2_1d} presents the variance evolution of $y_1$ and that of $y_2$ estimated by the quasi-Monte Carlo Sobol (MC-SOBOL) sequence with $1,000,000$ samples, AMR-GA with polynomial bases up to order 5 under various accuracy control values $TOL_1$ and global gPC Galerkin method with order $30$. Failure to capture the properties via the global gPC expansion after a short time $(t>=8)$ is also shown in Figure \ref{fig:KO_y1y2_1d}. Meanwhile, the results from AMR-GA converge. If we keep the order of the expansion fixed and make the tolerance $TOL_1$ stricter (smaller), more elements are desired via AMR-GA and higher accuracy is achieved.

\begin{table} [htbp]
\begin{center}
\begin{tabular}{llc|c|c|c}
\toprule
&AMR-CO&&CO&MC&MC-SOBOL\\
\midrule
        & No. of points & Error  &Error&Error&Error\\
\midrule
$TOL_1=10^{-3}$ & $160$ & $4.6e-2$ &$1.7e-1$ & $3.0e-1$ &$2.7e-1$\\
$TOL_1=10^{-4}$ & $260$ & $4.1e-3$ &$1.3e-1$& $2.4e-1$ &$1.4e-1$\\
$TOL_1=10^{-5}$ & $320$ & $2.8e-4$ &$1.0e-1$ & $9e-2$&$1.1e-1$\\
\bottomrule
\end{tabular}
\caption{Maximum relative errors of the variance when $t\in(0,30]$ for the Kraichnan-Orszag three-mode system with 1D random input via different method}
\label{tab:KO_1d_c}
\end{center}
\end{table}

\begin{table} [htbp]
\begin{center}
\begin{tabular}{l|lc|lc|lc}
\toprule
        & N & Error  &N &Error&N &Error\\
\midrule
$TOL_1$&$10^{-3}$&&$10^{-5}$&&$10^{-7}$\\
\midrule
$p=7$ & $30$ & $1.7e-1$ &$44$ & $1.8e-4$&$86$&$5.0e-6$ \\
$p=9$ & $20$ & $9.7e-2$ &$34$& $2.1e-4$&$62$&$6.8e-7$ \\
$p=11$ & $30$ & $1.1e-1$ &$34$ & $3.7e-4$&$52$&$1.7e-6$\\
\bottomrule
\end{tabular}
\caption{Maximum relative errors of the variance when $t\in(0,30]$ for the Kraichnan-Orszag three-mode system with 1D random input via AMR-GA}
\label{tab:KO_1d_a}
\end{center}
\end{table}

\begin{table} [htbp]
\begin{center}
\begin{tabular}{l|lc|lc|lc}
\toprule
        & N & Error  &N &Error&N &Error\\
\midrule
$TOL_1$&$10^{-3}$&&$10^{-5}$&&$10^{-7}$\\
\midrule
$p=7$ & $182$ & $3.8e-2$ &$352$ & $2.2e-4$ &$688$&$5.1e-6$\\
$p=9$ & $160$ & $4.6e-2$ &$320$& $2.8e-4$ &$640$&$9.9e-7$\\
$p=11$ & $216$ & $8.4e-2$ &$312$ & $2.4e-4$ &$624$&$2.0e-6$\\
\bottomrule
\end{tabular}
\caption{Maximum relative errors of the variance when $t\in(0,30]$ for the Kraichnan-Orszag three-mode system with 1D random input via AMR-CO}
\label{tab:KO_1d_b}
\end{center}
\end{table}

\begin{figure}[ht]
   \centering
   \subfigure[]{%
   \includegraphics[width = 5.5cm]{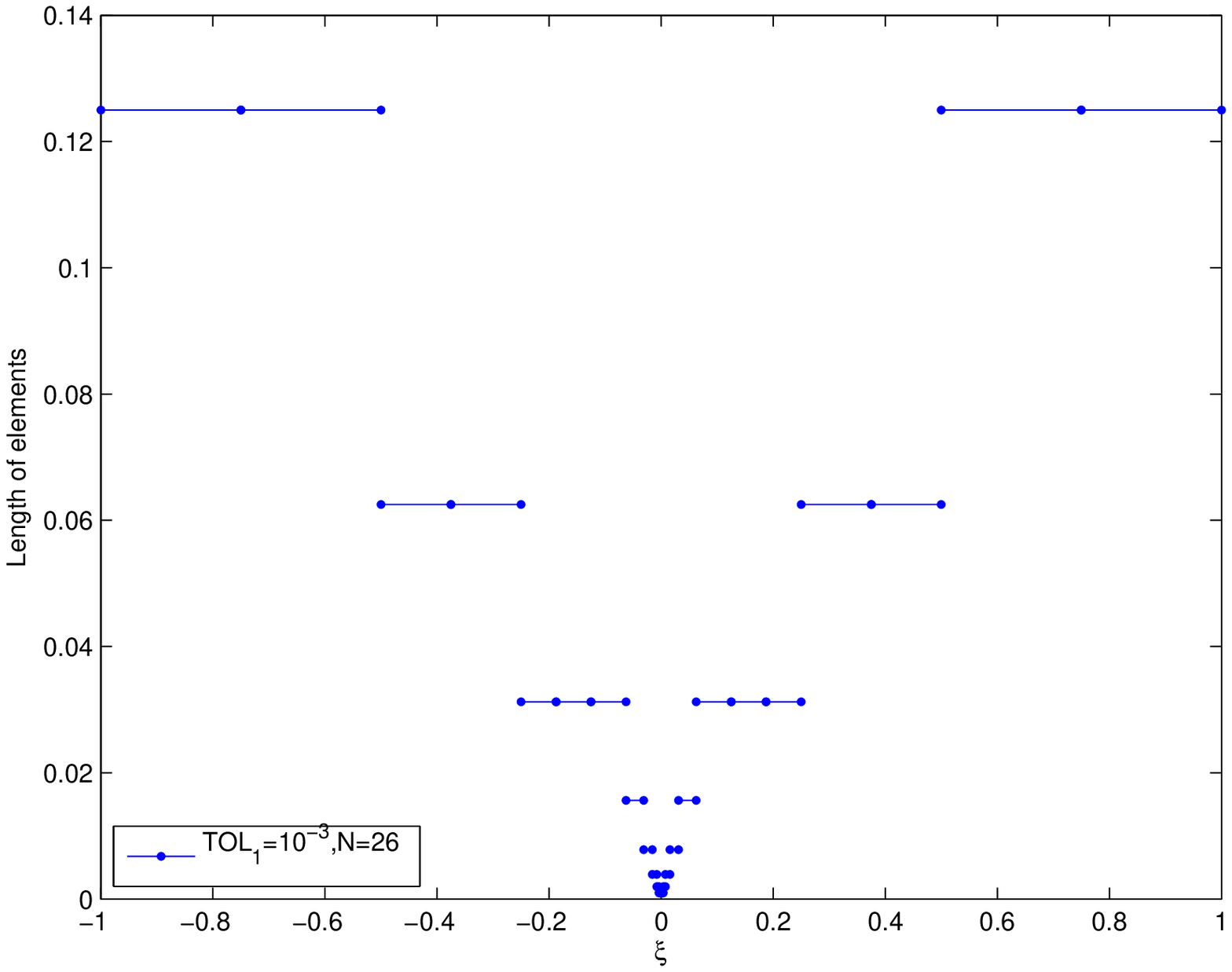}
   \label{KO_mesh_1d3_a}}
   \quad
   \subfigure[]{%
   \includegraphics[width = 5.5cm]{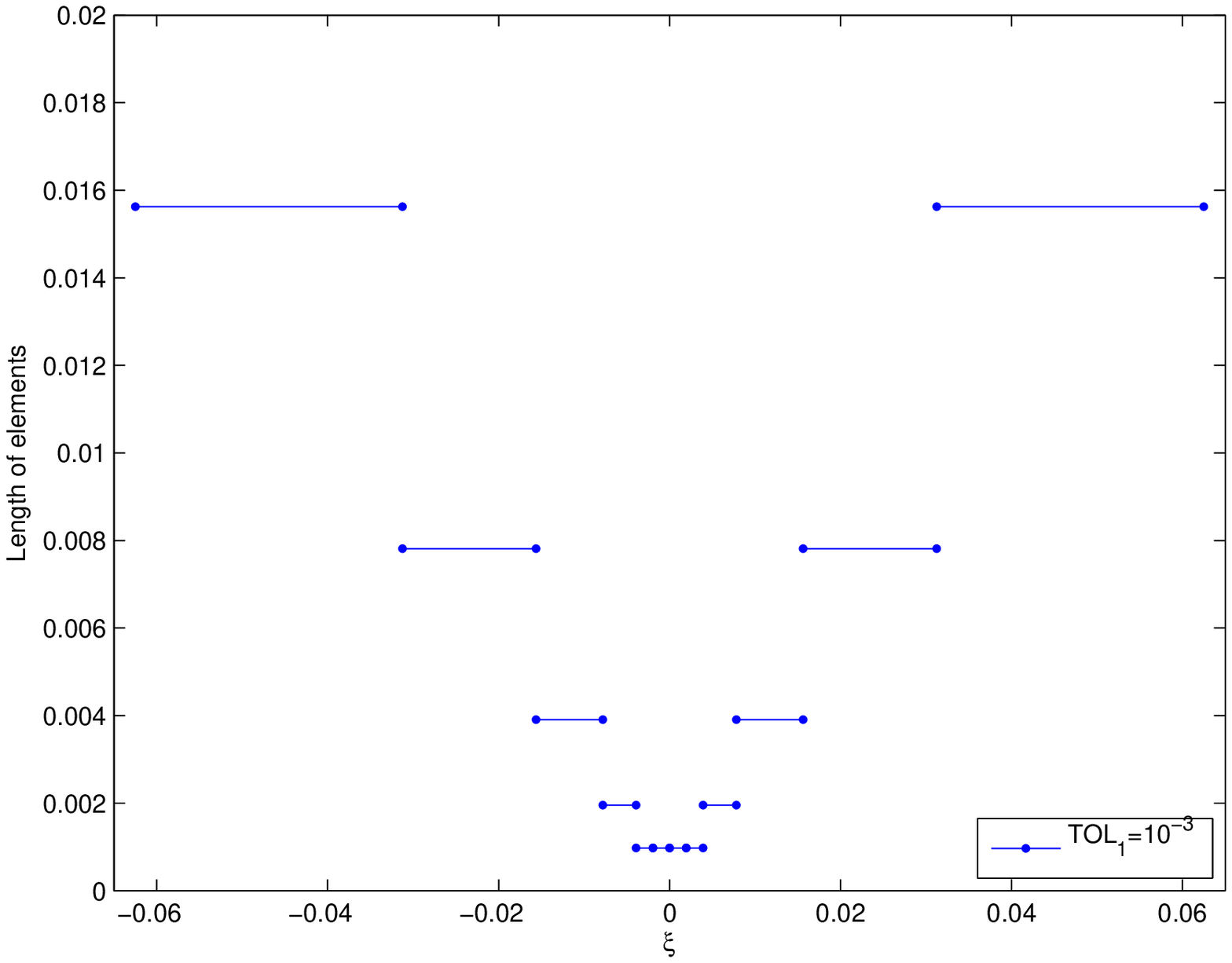}
   \label{KO_mesh_1d3_b}}
\caption{(a) Adaptive meshes for the Kraichnan-Orszag three-mode system with 1D random input when $p=7$ and $TOL_1=10^{-3}$ via AMR-GA; (b) zoom-in mesh of (a) near $\xi=0$.  }
\label{fig:KO_mesh_1d_ga}
\end{figure}

\begin{figure}[ht]
   \centering
   \subfigure[]{%
   \includegraphics[width = 5.5cm]{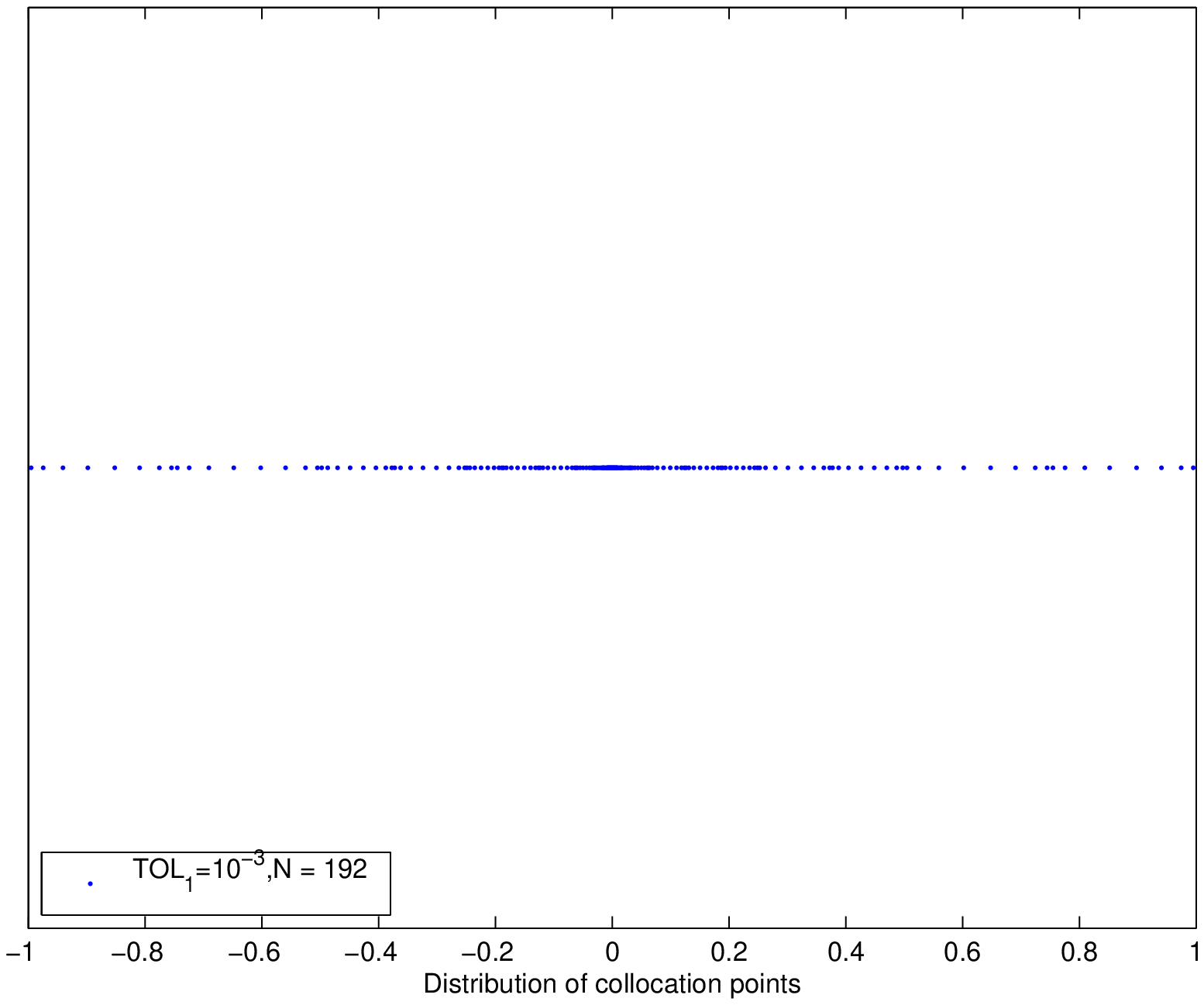}
   \label{KO_mesh_1d6_a}}
   \quad
   \subfigure[]{%
   \includegraphics[width = 5.5cm]{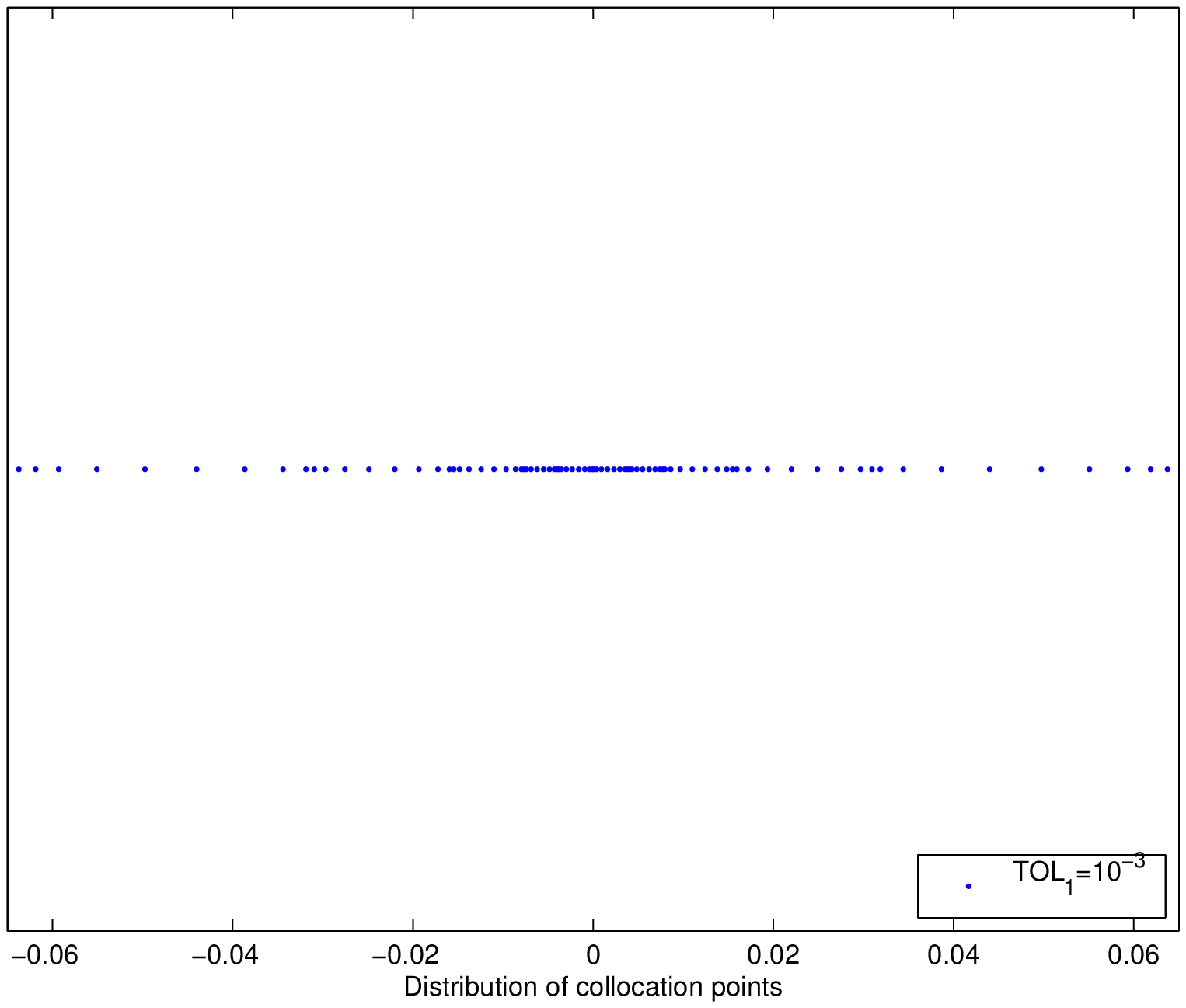}
   \label{KO_mesh_1d6_b}}
\caption{(a) Distribution of collocation points when solving the Kraichnan-Orszag three-mode system with 1D random input via AMR-CO with $p=7$ and $TOL_1=10^{-3}$; (b) zoom-in mesh of (a) near $\xi=0$.  }
\label{fig:KO_mesh_1d_co}
\end{figure}

We consider the solution given by MC-SOBOL with $1,000,000$ samples as the reference (ref) to study the error of AMR-CO, the global collocation method (CO), brute force Monte Carlo (MC) and MC-SOBOL. Results from AMR-CO with $p=9$ and different accuracy control values($TOL_1$) are presented in Table \ref{tab:KO_1d_c}. The error is defined as:
\begin{equation}\label{ex2_error}
\textrm{Error} =\sup_{t\in(0,30]}\max_{i=1,2,3} \frac{|{\rm{var}}(y_i)(t)-{\rm{var}}_{{\rm{ref}}}(y_i)(t)|}{{\rm{var}}_{{\rm{ref}}}(y_i)(t)}.
\end{equation}
In the first column of Table \ref{tab:KO_1d_c} we present the number of collocation points and the maximum of the relative errors of the variance of $y_1$, $y_2$ and $y_3$ when the system evolves to $t=30$ via AMR-CO with full model orders 9 under different $TOL_1$. In the implementation of the global collocation method we use Legendre collocation points with the same number of points as in AMR-CO. For MC and MC-SOBOL, the sample size equals the number of the collocation points in AMR-CO. Errors are shown in column two, three and four respectively. We note a much higher convergence rate for AMR-CO than for the other methods.

Secondly, since for this example the mesh refinement method appears to have a higher convergence rate we choose the result from AMR-CO with order 19 and very strict accuracy control $TOL_1 = 10^{-9}$ as the reference to study the relative error of both adaptive mesh refinement methods (AMR-GA and AMR-CO) with set of different orders and $TOL_1$. Table \ref{tab:KO_1d_a} shows the results of AMR-GA while Table \ref{tab:KO_1d_b} presents the results of AMR-CO. For a fixed value of the tolerance $TOL_1,$ higher order models require fewer elements in AMR-GA. However if we check carefully the number of variables we resolve in AMR-GA at $t = 30$ for each order, the magnitude of numbers are comparable for fixed $TOL_1$. For example, when $TOL_1 = 10^{-7}$, we solve $688$ variables at $t= 30$ for order 7 model, while that number becomes $620$ for order 9 model and $624$ for order 11 model respectively, and the relative errors of all three orders are at level $10^{-6}$. Furthermore, the same trend appears in the results of AMR-CO.

Details of the adaptive meshes resulting from our AMR-GA and AMR-CO algorithms around $\xi=0$ are presented in Figures \ref{fig:KO_mesh_1d_ga} and \ref{fig:KO_mesh_1d_co}. The finest meshes are around the discontinuity of the random space. This shows that our mesh refinement criterion locates accurately the discontinuity. Furthermore, when $TOL_1$ is extremely small, the meshes exhibit the pattern that the closer the element is to $\xi=0,$ the smaller the element in AMR-GA and the more the collocation points in AMR-CO. Meanwhile the meshes and points are symmetric with respect to $\xi=0$ as they should be according to the symmetry of the K-O system.

\subsubsection{Two-dimensional random input}
\begin{figure}[ht]
   \centering
   \subfigure[]{%
   \includegraphics[width = 5.5cm]{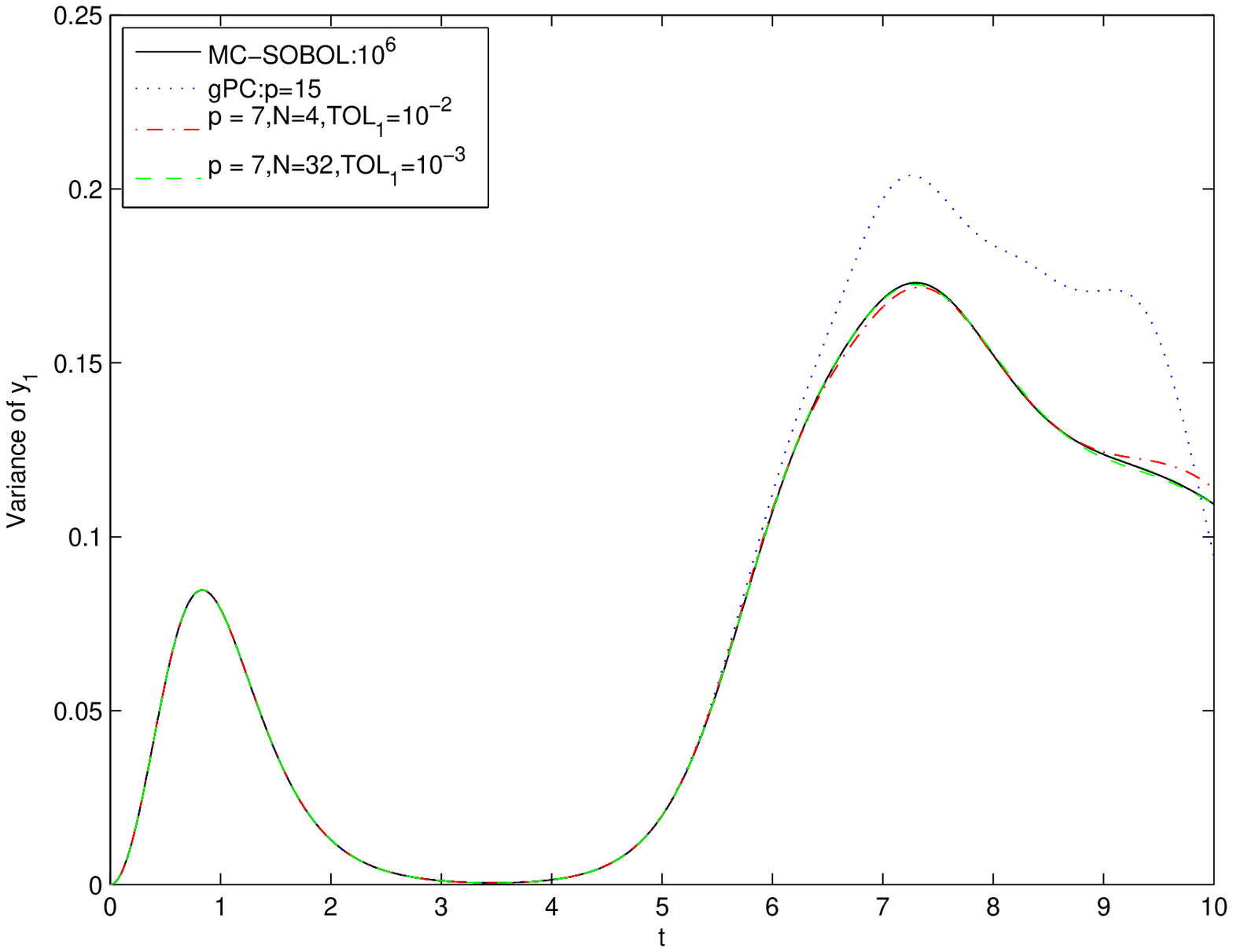}
   \label{y1_var_2d}}
   \subfigure[]{%
   \includegraphics[width = 5.5cm]{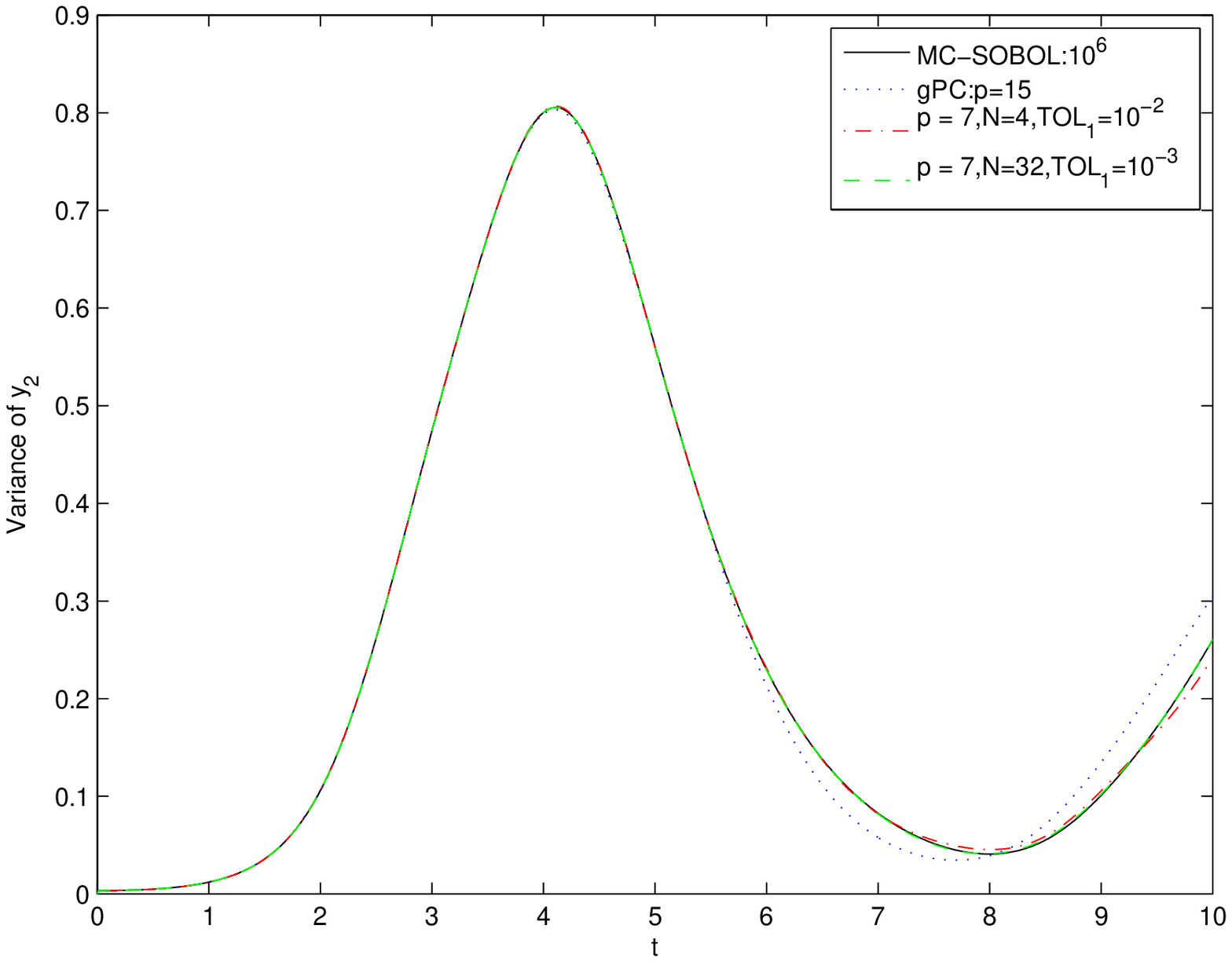}
   \label{y2_var_2d}}
   \subfigure[]{%
   \includegraphics[width = 5.5cm]{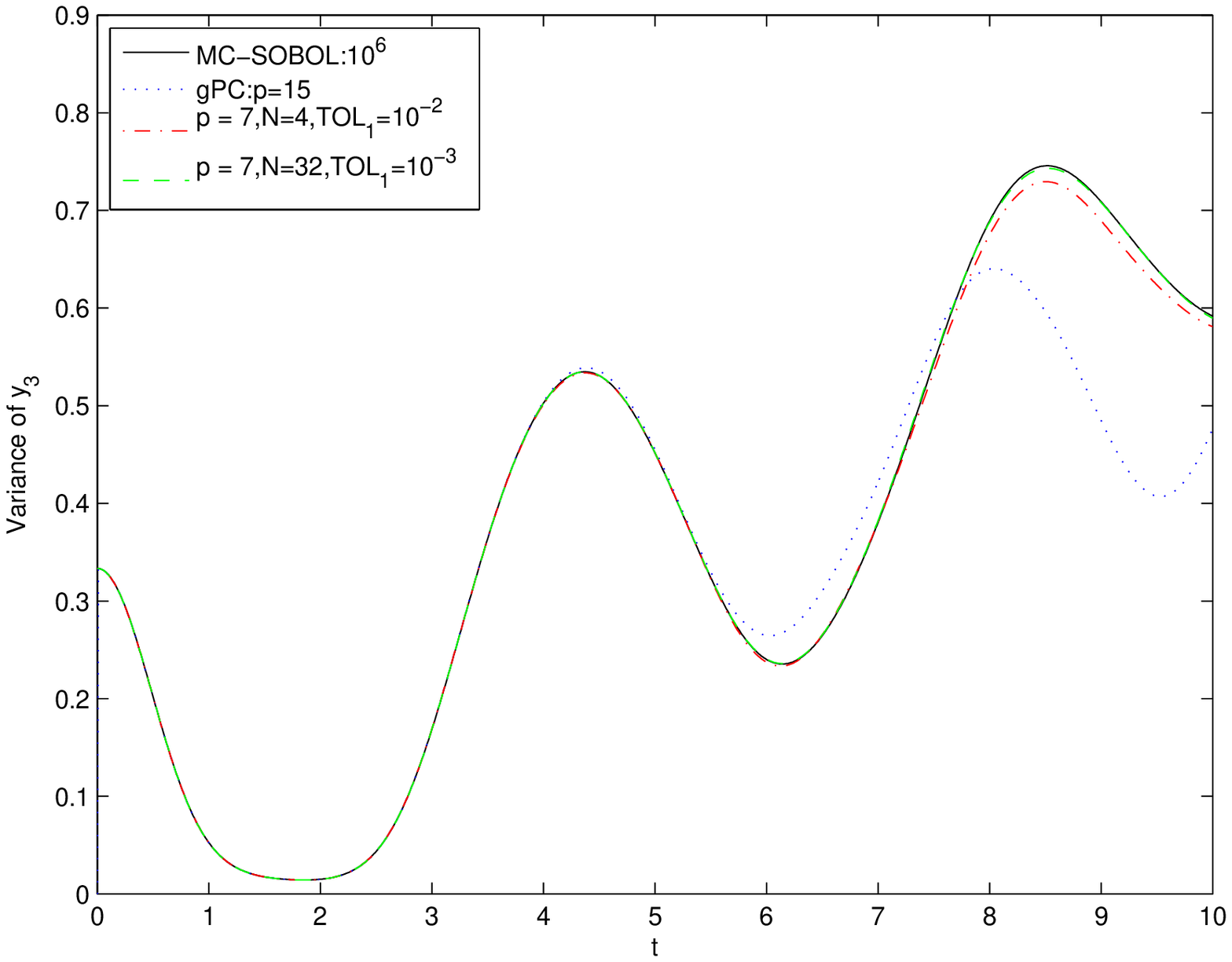}
   \label{y3_var_2d}}
   \subfigure[]{%
   \includegraphics[width = 5.5cm]{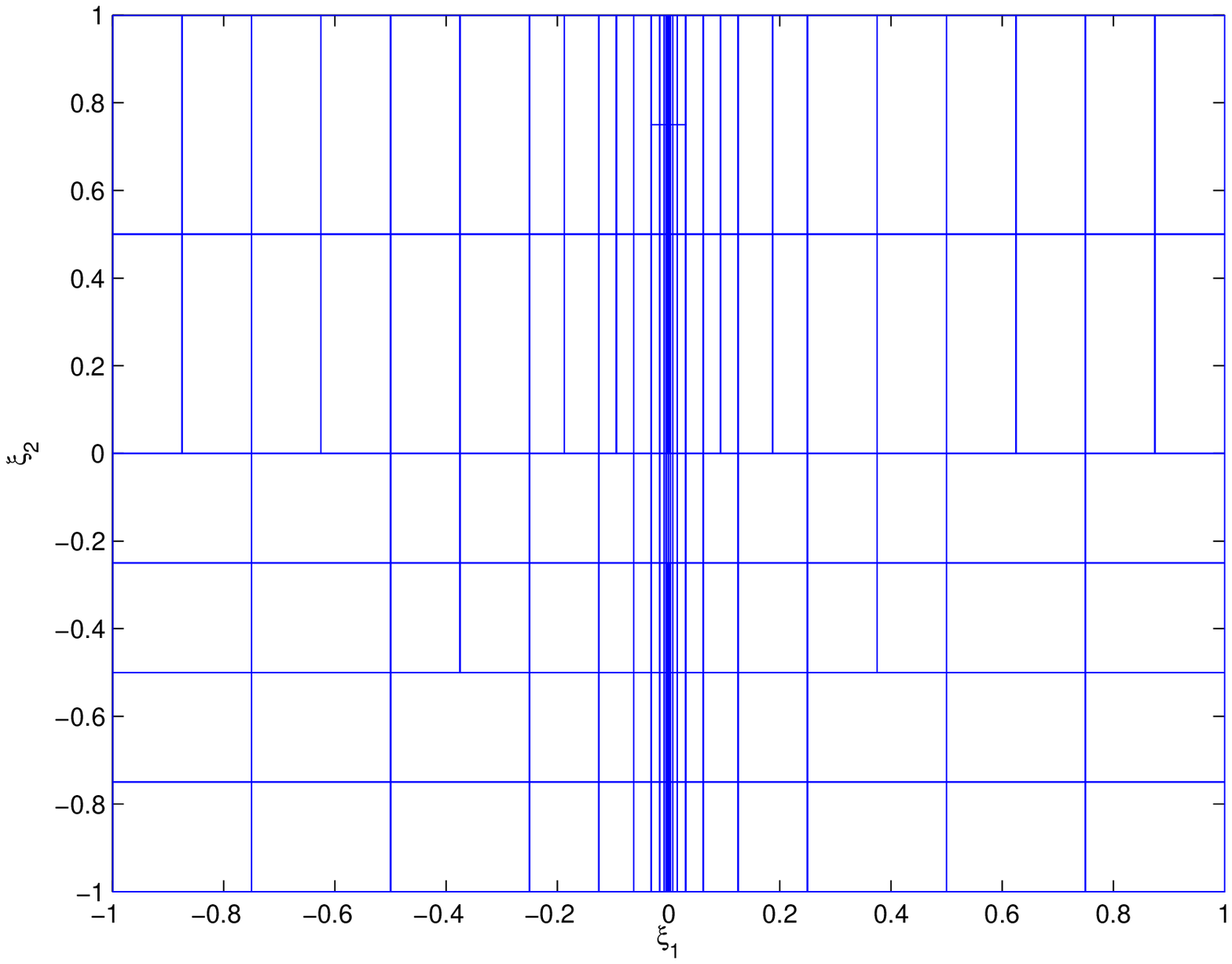}
   \label{KO_mesh_2d}}
\caption{The Kraichnan-Orszag three-mode system with 2D random inputs, $p=7,TOL_2=0.1$. (a) Evolution of variance of $y_1$; (b) Evolution of variance of $y_2$; (c) Evolution of variance of $y_3$; (d)Adaptive meshes for 2D random input when $p=7, TOL_1=10^{-6}, N=210$.
  }
\label{fig:KO_mesh_2d}
\end{figure}

We study the K-O system with initial conditions involving two independent random inputs
\begin{equation}\label{ex:KO_ini_2d}
y_1(0;\o) = 1,\quad y_2(0;\o) = 0.1\xi_1(\o),\quad y_3(0;\omega)=\xi_2(\o),
\end{equation}
where $\xi_1$ and $\xi_2$ are independent uniform random variables on $[-1,1]$.
In Figure \ref{fig:KO_mesh_2d}, we plot the evolution of the variance of each random output $y_i,i=1,2,3$ subject to a 2D random input obtained by AMR-GA with order 7 for the full model with different values of $TOL_1.$ We also show the mesh of the random space at time $t=10$ for $210$ elements generated by AMR-GA with order 7 for the full model and, $TOL_1=10^{-6}$ and $TOL_2=0.1$. The smallest elements are around the discontinuity $y_2=0$ and the results are more sensitive to $\xi_1$ because of the discontinuity introduced by $\xi_1$. 

To study the convergence of AMR-CO compared with CO, MC and MC-SOBOL, we take the result from MC-SOBOL with $1,000,000$ samples as the reference. To implement the direct collocation method, we choose the tensor grid with $\lceil {N}^{\frac{1}{2}} \rceil$ if $\lceil {N}^{\frac{1}{2}} \rceil$ is even, otherwise we choose $\lceil {N}^{\frac{1}{2}} \rceil+1$ Legendre quadrature points in each dimension. As for MC and MC-SOBOL, we take the sample size to be the number of collocation points in AMR-CO. Table \ref{tab:KO_2d_a} demonstrates the relative error with respect to the reference. The results show that AMR-CO converges much faster than the other three methods.

The result from AMR-CO with order $10$ for the full model, $TOL_1 = 10^{-9}$ and $TOL_2=0.1$ is selected to be the reference to derive the relative errors of low ordered models. Table \ref{tab:KO_2d_b} presents the relative errors of variance by AMR-GA of different sets of full model orders and different levels of accuracy control, while Table \ref{tab:KO_2d_c} presents the relative errors of variance by AMR-CO. We observe similar trends as in the 1D case, namely that more accurate models require fewer elements in AMR-GA if the accuracy tolerance is kept fixed. Also, that stricter accuracy control requires more elements if the order of the models is fixed. To achieve the same level of relative errors, the number of the elements increases faster in the 2D case than the 1D case. However, in AMR-CO with the same accuracy control, the number of required collocation points are in the same order of magnitude for different full model orders. Also, the relative errors are in the same order of magnitude. 
\begin{table} [htbp]
\begin{center}
\begin{tabular}{llc|c|c|c}
\toprule
&AMR-CO&&CO&MC&MC-SOBOL\\
\midrule
        & No. of points & Error  &Error&Error&Error\\
\midrule
$TOL_1=10^{-2}$ & $576$ & $1.8e-1$ &$3.7e-1$ & $1.5e-1$ &$2.5e-1$\\
$TOL_1=10^{-3}$ & $1944$ & $3.7e-2$ &$1.7e-1$& $1.0e-1$ &$2.1e-1$\\
$TOL_1=10^{-4}$ & $4896$ & $3.1e-3$ &$8.5e-2$ & $6.2e-2$&$9.3e-2$\\
\bottomrule
\end{tabular}
\caption{Maximum relative errors of the variance when $t\in(0,10]$ for the Kraichnan-Orszag three-mode system with 2D random input via different method}
\label{tab:KO_2d_a}
\end{center}
\end{table}

\begin{table} [htbp]
\begin{center}
\begin{tabular}{l|lc|lc|lc}
\toprule
        & N & Error  &N &Error&N &Error\\
\midrule
$TOL_1$&$10^{-3}$&&$10^{-5}$&&$10^{-7}$\\
\midrule
$p=5$ & $34$ & $2.8e-2$ &$222$ & $2.6e-3$&$424$&$1.7e-4$ \\
$p=7$ & $32$ & $2.7e-2$ &$152$ & $7.7e-4$&$310$&$4.7e-6$\\
\bottomrule
\end{tabular}
\caption{Maximum relative errors of the variance when $t\in(0,10]$ for the Kraichnan-Orszag three-mode system with 2D random input via Galerkin approach}
\label{tab:KO_2d_b}
\end{center}
\end{table}

\begin{table} [htbp]
\begin{center}
\begin{tabular}{l|lc|lc|lc}
\toprule
        & N & Error  &N &Error&N &Error\\
\midrule
$TOL_1$&$10^{-3}$&&$10^{-5}$&&$10^{-7}$\\
\midrule
$p=5$ & $1944$ & $3.7e-2$ &$10224$ & $3.2e-4$&$26496$&$3.8e-6$ \\
$p=7$ & $2048$ & $7.4e-2$ &$9728$ & $8.9e-4$&$19840$&$6.4e-6$\\
\bottomrule
\end{tabular}
\caption{Maximum relative errors of the variance when $t\in(0,10]$ for the Kraichnan-Orszag three-mode system with 2D random input via collocation approach}
\label{tab:KO_2d_c}
\end{center}
\end{table}

\subsubsection{Three-dimensional random input}
\begin{figure}[htbp]
   \centerline{
   \psfig{file=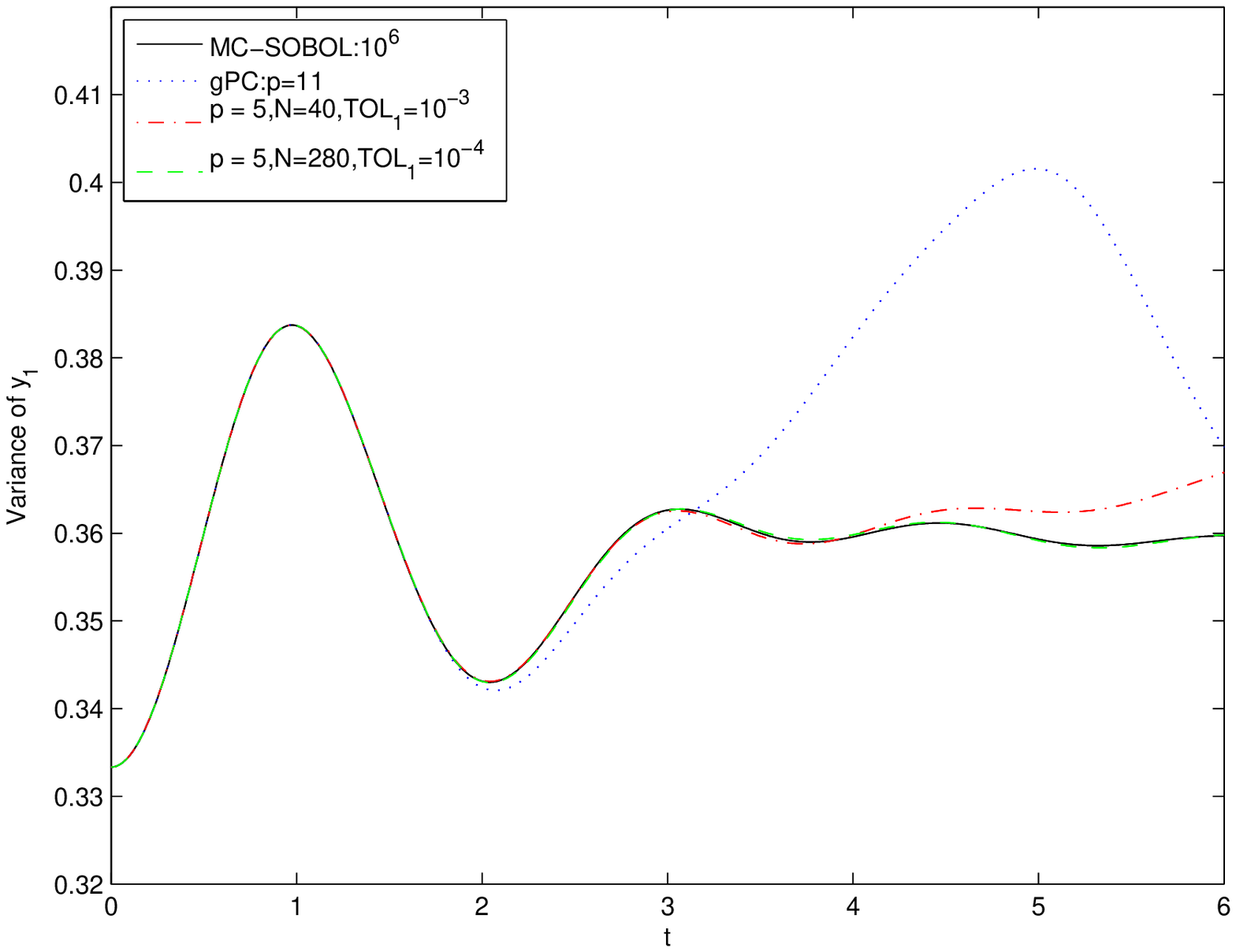,width=7cm}
   \psfig{file=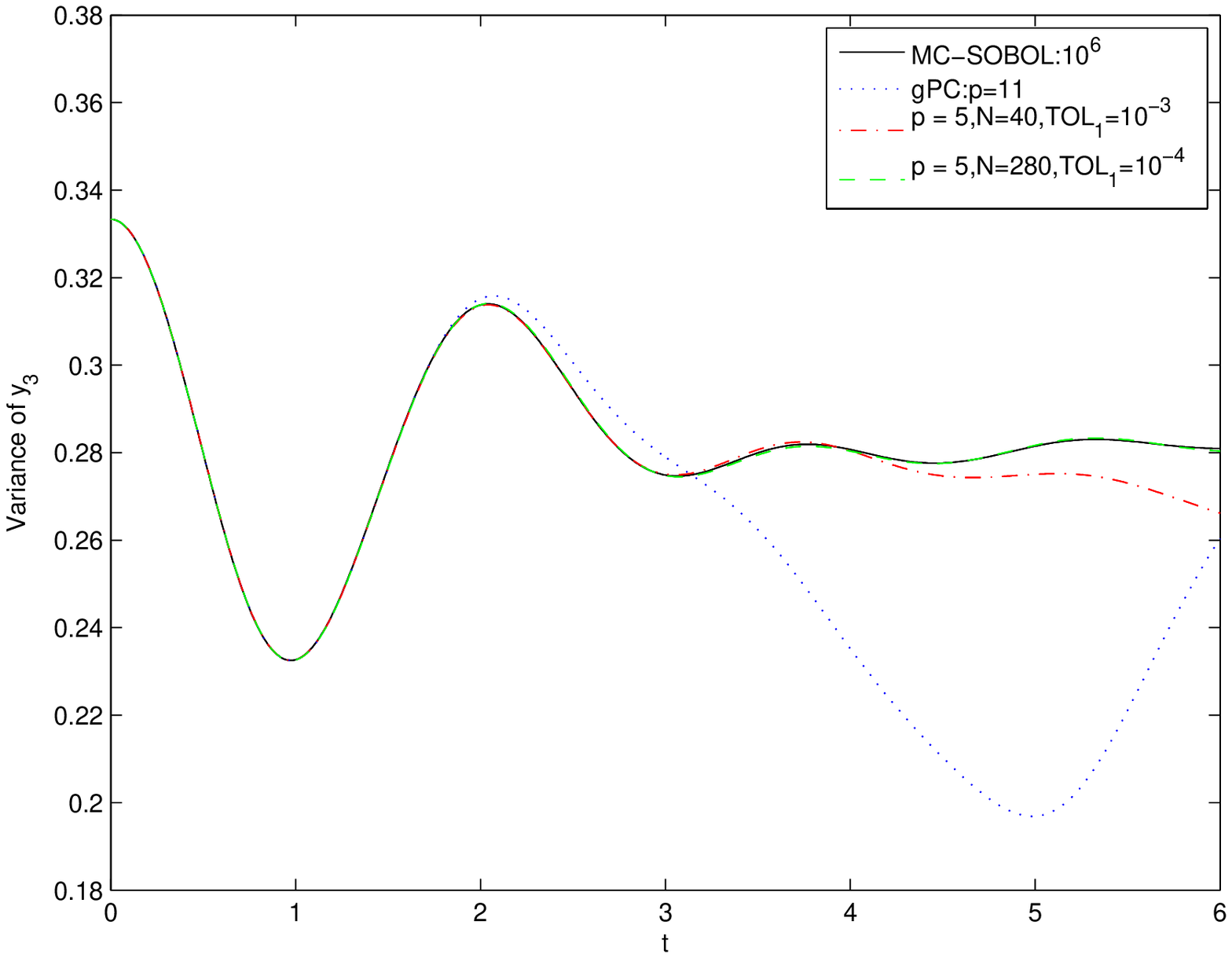,width=7cm}
  }
\caption{Evolution of the variance of $y_1$(left) and evolution of the variance of $y_3$(right) for the Kraichnan-Orszag three-mode system with 3D random inputs.
  }
\label{fig:KO_mesh_3d}
\end{figure}

The initial conditions in this case are
\begin{equation}\label{ex:KO_ini_3d}
y_1(0;\o) = \xi_1(\o), \quad y_2(0;\o) = \xi_2(\o),\quad y_3(0;\omega)=\xi_3(\o),
\end{equation}
where $\xi_1$, $\xi_2$ and $\xi_3$ are independent uniform random variables on $[-1,1]$. In this case, discontinuities occur at $y_1=0$ and $y_2=0.$ Due to the presence of discontinuities and moderately high dimensionality, this problem is more difficult than the previous ones. Figure \ref{fig:KO_mesh_3d} shows the evolution of the variance of $y_1$ and $y_3$ obtained by an order 5 full model with different $TOL_1$. The variances of $y_1$ and $y_2$ are the same due to the symmetry of $y_1$ and $y_2$ in \eqref{ex:KO} and the corresponding initial random input. The results for a global gPC expansion of order 9 diverges from the MC-SOBOL results at $t\approx 3.$ On the other had, our AMR-GA algorithms achives much better accuracy. Similarly as in the previous two cases, a comparison of the relative errors of the output variances between AMR-CO, global collocation, MC and MC-SOBOL is presented in Table \ref{tab:KO_3d_a}. (the MC-SOBOL result with $1,000,000$ samples is taken as the reference). Again, the mesh refinement algorithm shows a better convergence than the other methods, however the advantage is not as big as before. 

We choose results from AMR-CO with $p=10,TOL_1=10^{-7}$ as the reference and show the results of AMR-GA with different set of orders and $TOL_1$ in Table \ref{tab:KO_3d_b} and that of AMR-CO in Table \ref{tab:KO_3d_c}. As we can see, the number of the elements in AMR-GA as well as the collocation points in AMR-CO grows dramatically fast in order to gain sufficient accuracy compared to the 1D and 2D cases. Therefore, both of our adaptive mesh refinement algorithms are still dimension-dependent algorithms and many more elements or collocations points are required to resolve the discontinuity.

\begin{table} [htbp]
\begin{center}
\begin{tabular}{llc|c|c|c}
\toprule
&AMR-CO&&CO&MC&MC-SOBOL\\
\midrule
        & No. of points & Error  &Error&Error&Error\\
\midrule
$TOL_1=10^{-2}$ & $2784$ & $3.4e-2$ &$8.1e-2$ & $4.7e-2$ &$6.4e-2$\\
$TOL_1=10^{-3}$ & $4176$ & $2.4e-2$ &$7.0e-2$& $3.7e-2$ &$8.5e-2$\\
$TOL_1=10^{-4}$ & $23664$ & $3.4e-3$ &$3.8e-2$ & $1.6e-2$&$2.8e-2$\\
\bottomrule
\end{tabular}
\caption{Maximum relative errors of the variance when $t\in(0,6]$ for the Kraichnan-Orszag three-mode system with 3D random input via different method}
\label{tab:KO_3d_a}
\end{center}
\end{table}

\begin{table} [htbp]
\begin{center}
\begin{tabular}{l|lc|lc|lc}
\toprule
        & N & Error  &N &Error&N &Error\\
\midrule
$TOL_1$&$10^{-2}$&&$10^{-3}$&&$10^{-4}$\\
\midrule
$p=4$ & $32$ & $8.6e-2$ &$48$ & $3.0e-2$&$312$&$2.7e-3$ \\
$p=5$ & $32$ & $7.1e-2$ &$40$& $5.3e-2$&$280$&$3.0e-3$ \\
$p=6$ & $16$ & $9.9e-2$ &$32$ & $2.9e-2$&$112$&$1.7e-3$\\
\bottomrule
\end{tabular}
\caption{Maximum relative errors of the variance when $t\in(0,6]$ for the Kraichnan-Orszag three-mode system with 3D random input via Galerkin approach}
\label{tab:KO_3d_b}
\end{center}
\end{table}

\begin{table} [htbp]
\begin{center}
\begin{tabular}{l|lc|lc|lc}
\toprule
        & N & Error  &N &Error&N &Error\\
\midrule
$TOL_1$&$10^{-3}$&&$10^{-5}$&&$10^{-7}$\\
\midrule
$p=4$ & $2784$ & $3.4e-2$ &$4176$ & $2.4e-2$&$23664$&$3.4e-3$ \\
$p=5$ & $2160$ & $5.9e-2$ &$5400$& $7.9e-3$&$29160$&$3.1e-3$ \\
$p=6$ & $3312$ & $1.3e-1$ &$6624$ & $2.1e-2$&$21528$&$3.2e-3$\\
\bottomrule
\end{tabular}
\caption{Maximum relative errors of the variance when $t\in(0,6]$ for the Kraichnan-Orszag three-mode system with 3D random input via collocation approach}
\label{tab:KO_3d_c}
\end{center}
\end{table}

\subsection{Kuramoto-Sivashinsky equation}
In this section, we study the normalized Kuramoto-Sivashinsky (K-S) equation with periodic boundary conditions in the interval $[0,2\pi].$ Let the damping parameter to the original value derived by Sivashinsky, $\nu=4$, and introduce the bifurcation parameter $\alpha = 4\tilde{L}^2 = \frac{L^2}{4\pi}$. The equation can be restated as
\begin{equation}\label{k_s}
\begin{split}
& u_t=-4u_{xxxx}-\alpha\left[u_{xx}+\frac{1}{2}(u_x)^2\right],\quad 0\leq x\leq 2\pi,\\
& u(x+2\pi,t) = u(x,t),\quad u(x,0) = u_0(x).
\end{split}
\end{equation}
The K-S equation solution exhibits a broad range of behaviors as $\alpha$ varies. In particular, as $\alpha$ is varied the solution can exhibit bifurcations, thus a small change in $\alpha$ may result in a dramatic change of the solution. We examine the mean of the solution when the parameter $\alpha$ is in an interval which includes several bifurcation values. High precision is required for the numerical simulation due to the extreme sensitivity of the solution to the stepsize used \cite{hyman1986order,hyman1986kuramoto}. We use high precision pseudospectral approximations for the spatial derivatives, and semi implicit time differencing with fixed time step for the temporal derivative. We only present results for the adaptive mesh refinement algorithm with collocation (AMR-CO).

\begin{figure}[ht]
   \centering
   \subfigure[]{%
   \includegraphics[width = 5.5cm]{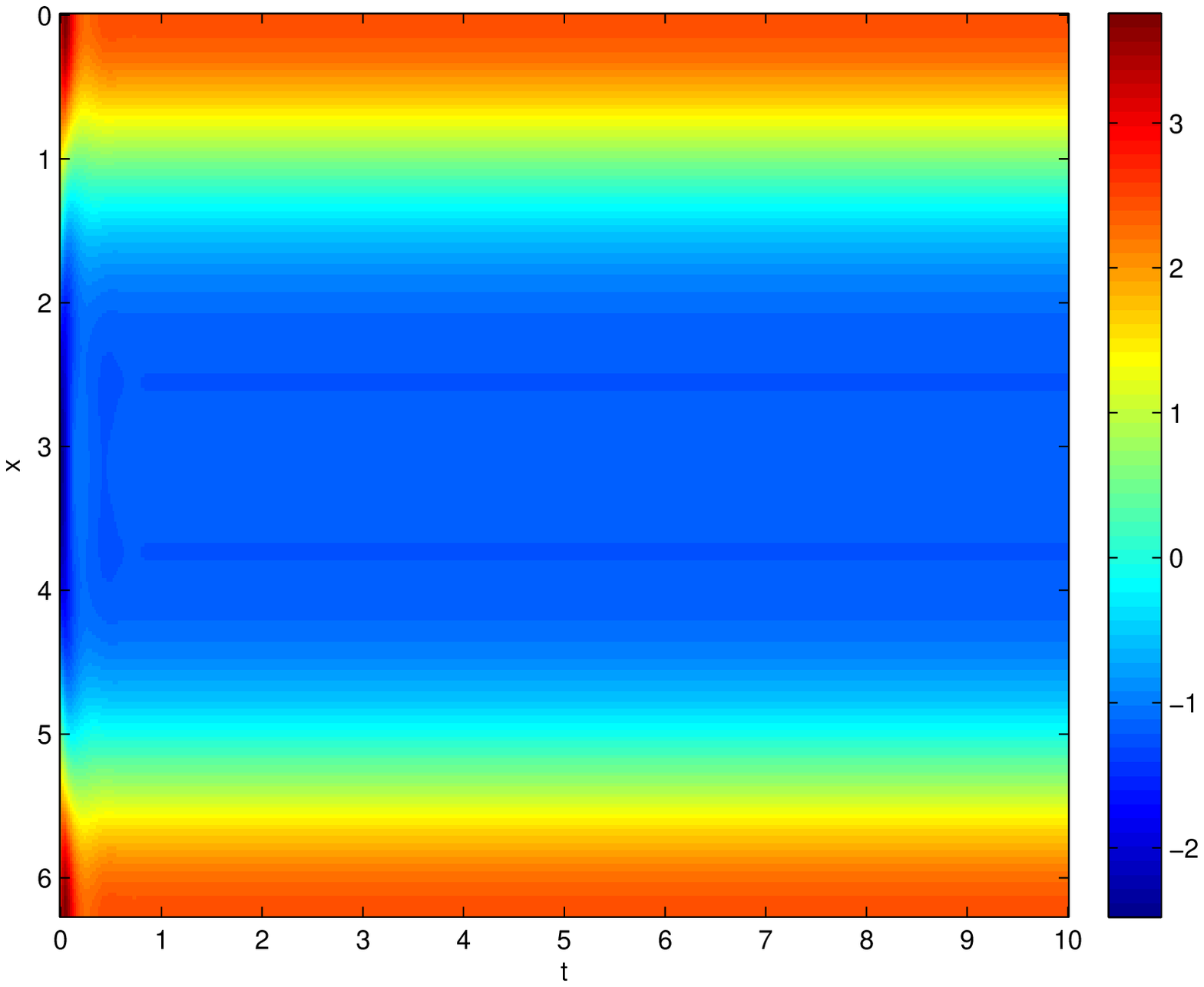}
   \label{fig:ks_13}}
   \subfigure[]{%
   \includegraphics[width = 5.5cm]{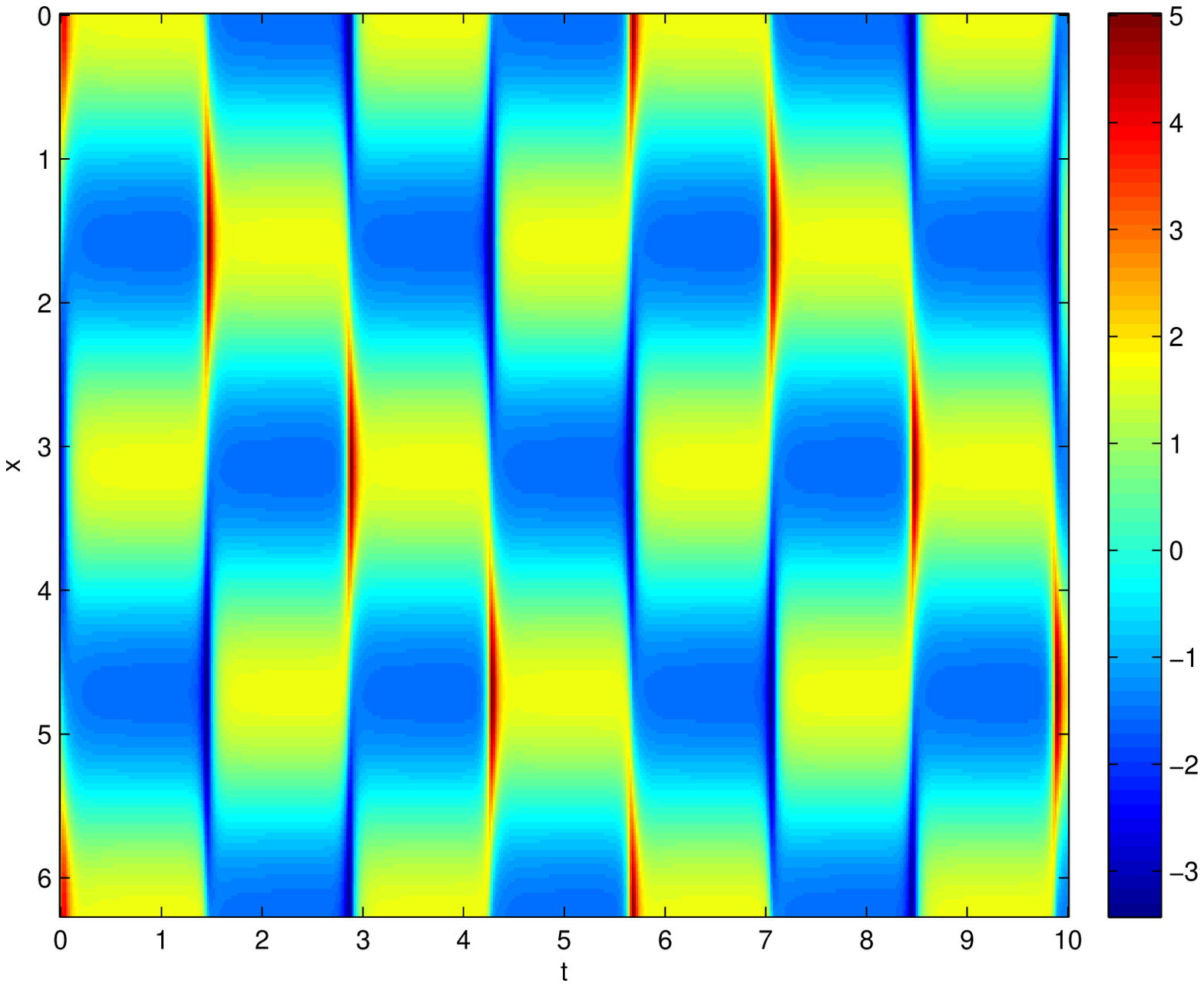}
   \label{fig:ks_17}}
\caption{Solutions to K-S equation with bifurcation parameter $\alpha = 13$(left) and $\alpha = 17$(right) when $t\in[0,10]$ }
\label{fig:ex3_1}
\end{figure}

We restrict attention to the interval $\alpha\in[13,17]$ in which the solution can be unimodal or bimodal depending on the value of $\alpha.$ We obtain the solution for the initial condition $u_0(x)  = 2.6680\cos(x)+0.1979\cos(2x)+0.0094\cos(3x)$. For $\alpha$ in the unimodal regime, the solution quickly converges to the unimodal steady solution. Figure \ref{fig:ks_13} shows the unimodal solution when $\alpha = 13$ while Figure \ref{fig:ks_17} presents the periodical bimodal solution when $\alpha = 17$. The solution to K-S equation is very sensitive with respect to the parameter $\alpha$, and we need to start with an adequate number of elements in random space. We start with $32$ initial elements in each refinement implementation. Since the expression of the RHS in \eqref{k_s} is complicated, only the decoupled AMR-CO is simulated. Result from MC-SOBOL simulation with $1,000,000$ samples is considered as the reference. 

\begin{figure}[ht]
   \centering
   \subfigure[]{%
   \includegraphics[width = 5.5cm]{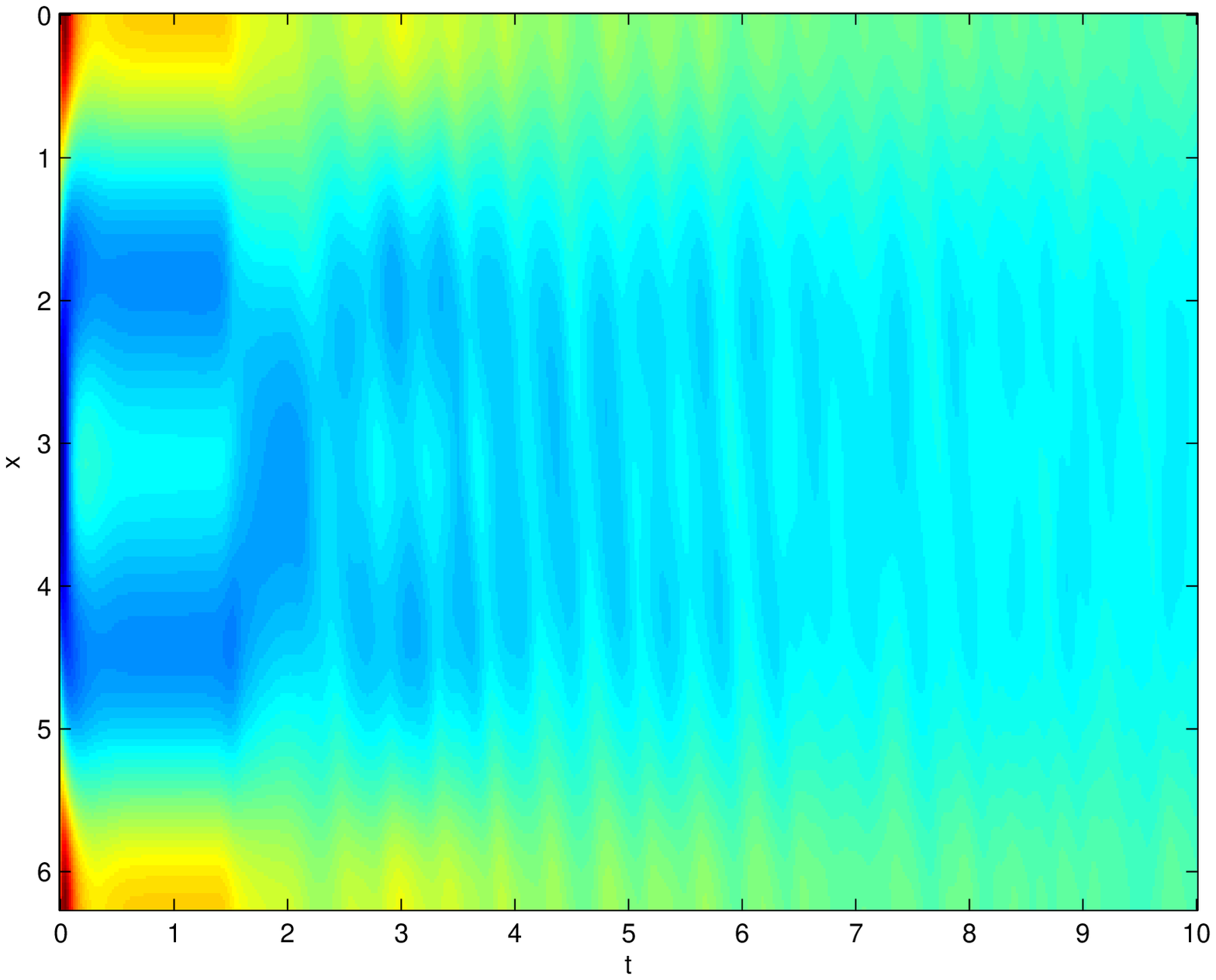}
   \label{fig:ks_mean_co}}
   \subfigure[]{%
   \includegraphics[width = 5.5cm]{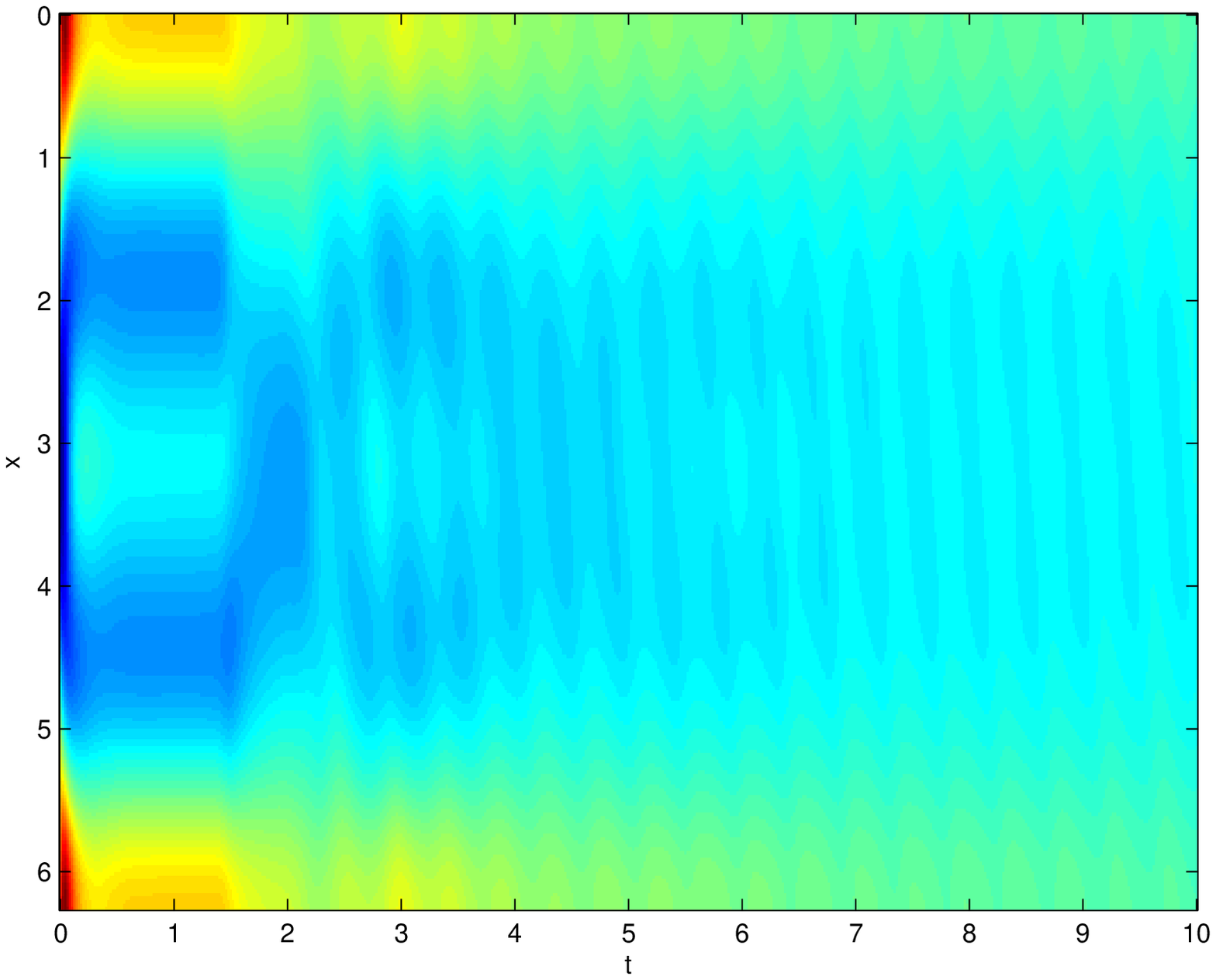}
   \label{fig:ks_mean_ref}}
\caption{Mean of solution to K-S equation with bifurcation parameter $\alpha\sim U[13,17]$ when $t\in[0,10]$ via AMR-CO with $d = 11$, $TOL_1=0.1$(left) and MC-SOBOL with $1,000,000$ samples(right)}
\label{fig:ex3_2}
\end{figure}

\begin{figure}[ht]
   \centering
   \subfigure[]{%
   \includegraphics[width = 5.5cm]{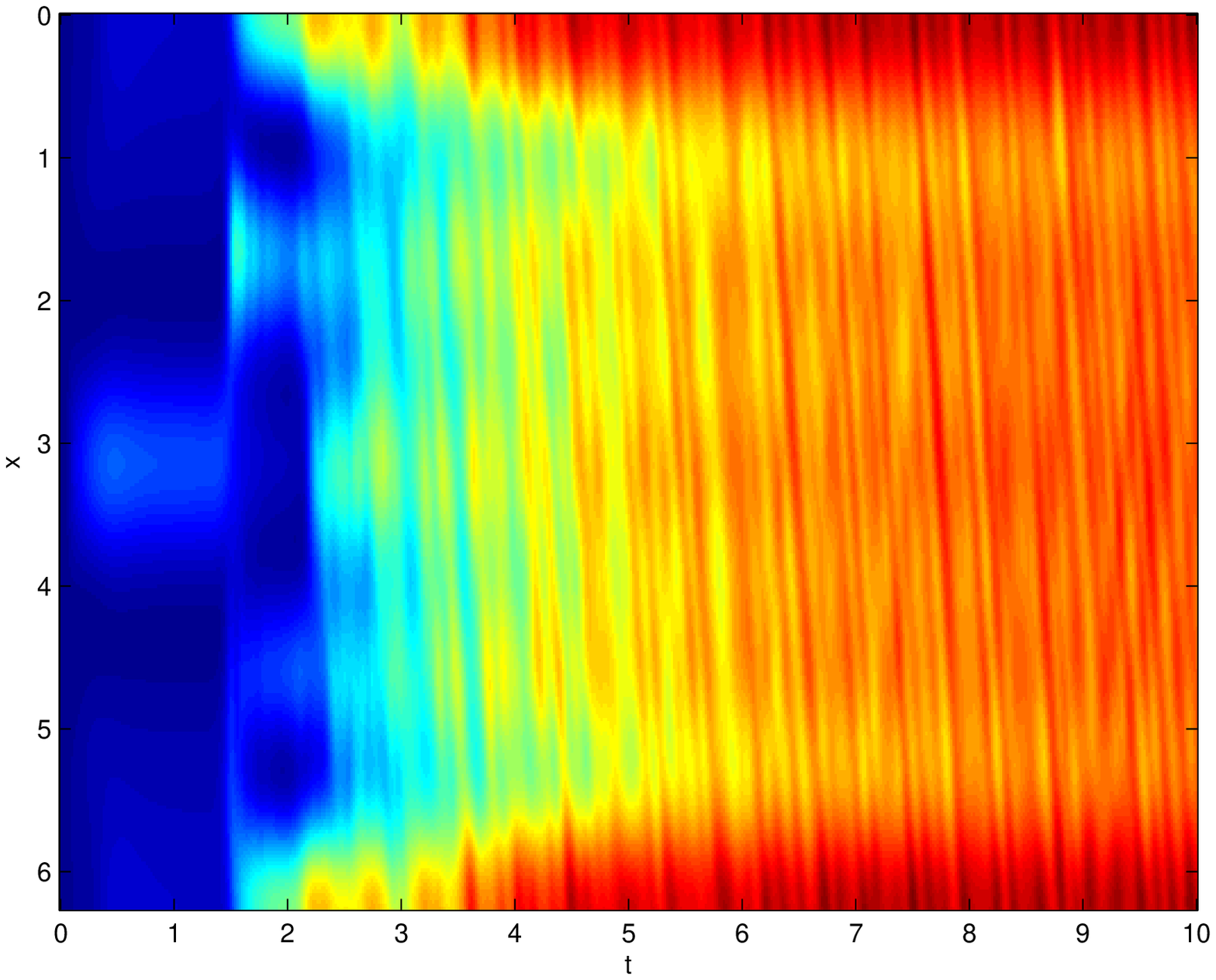}
   \label{fig:ks_var_co}}
   \subfigure[]{%
   \includegraphics[width = 5.5cm]{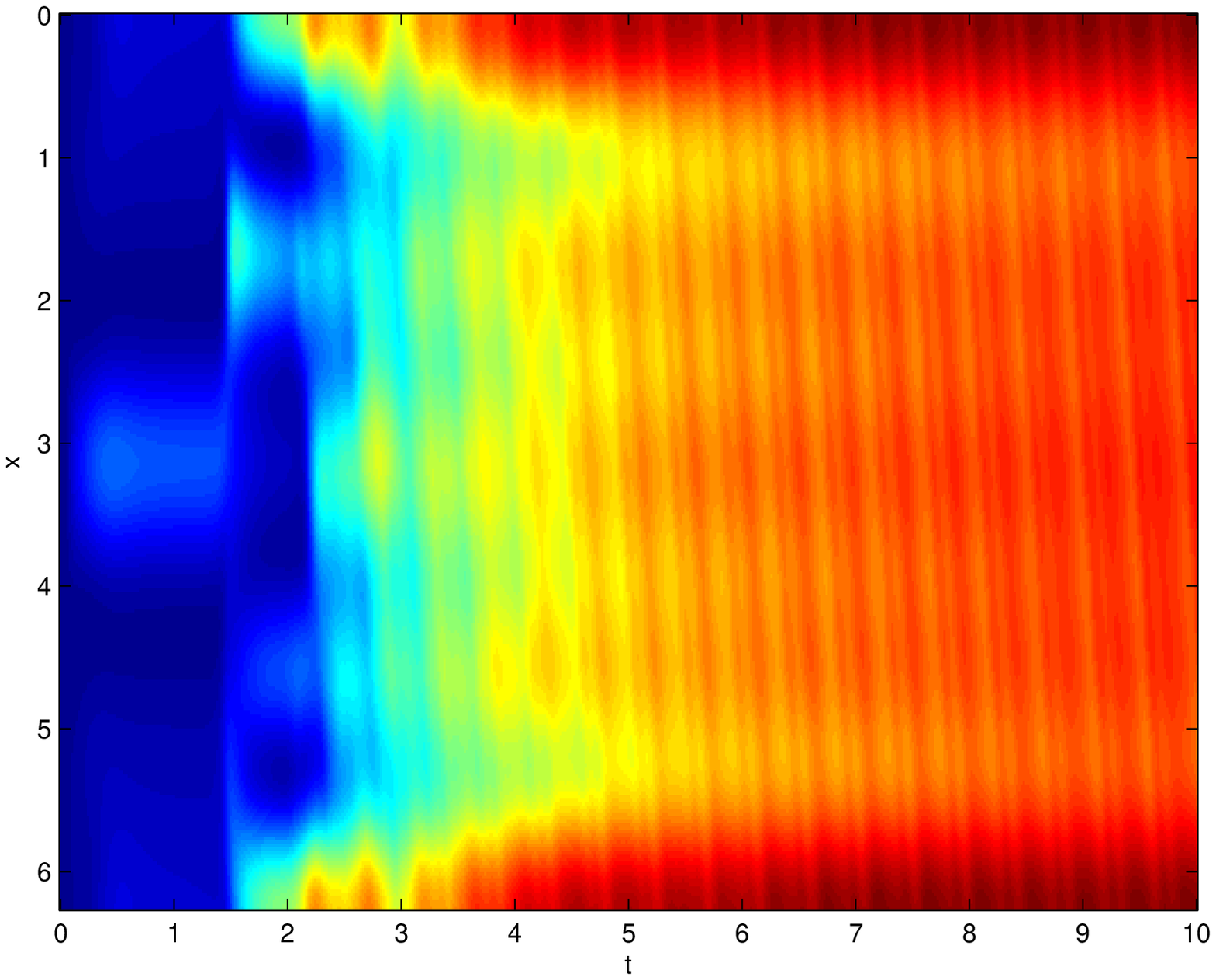}
   \label{fig:ks_var_ref}}
\caption{Variance of solution to K-S equation with bifurcation parameter $\alpha\sim U[13,17]$ when $t\in[0,10]$ via AMR-CO with $d = 11$, $TOL_1=0.1$(left) and MC-SOBOL with $1,000,000$ samples(right)}
\label{fig:ex3_3}
\end{figure}

\begin{figure}[htbp]
   \centering
   \subfigure[]{%
   \includegraphics[width=5.5cm]{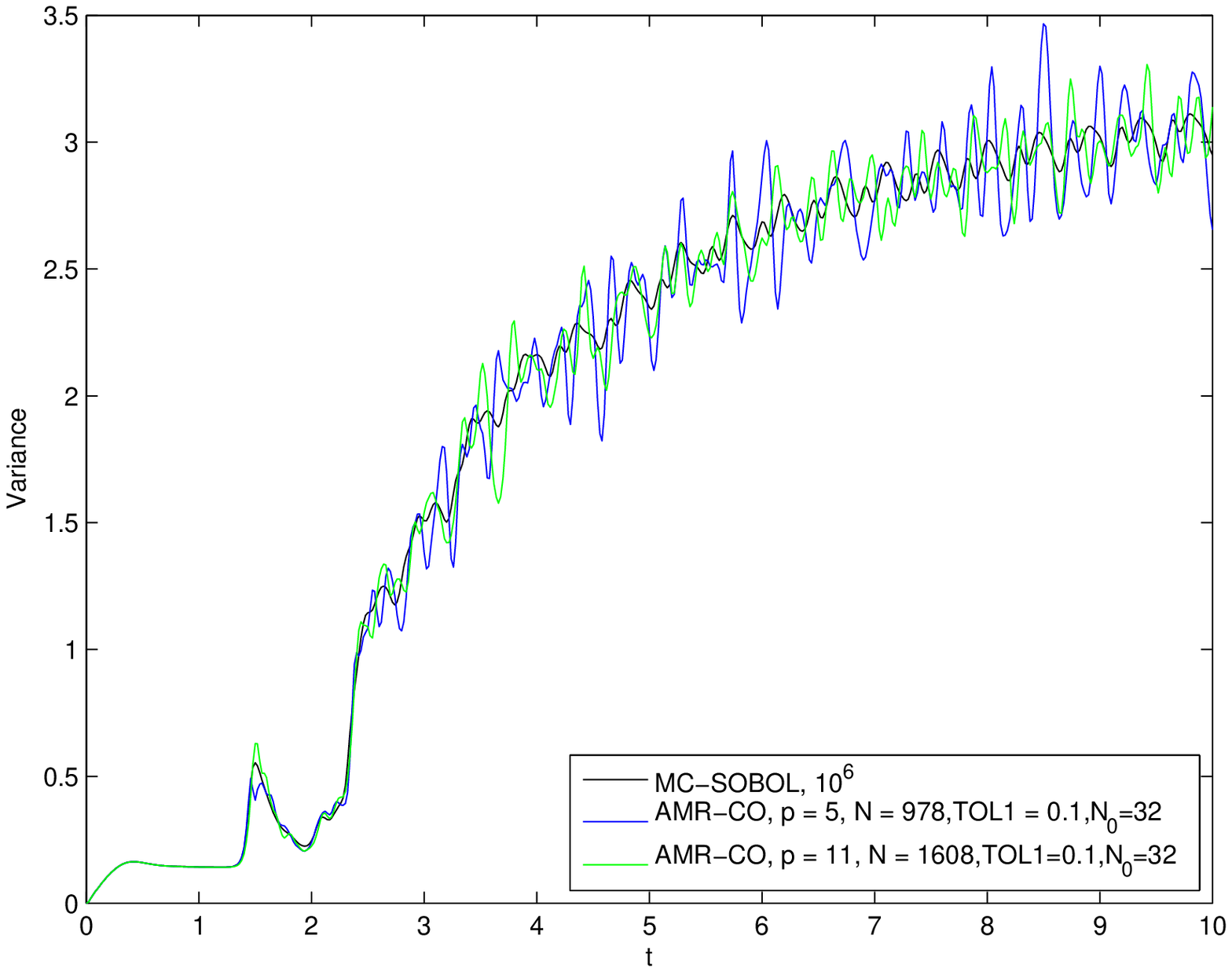}
   \label{fig:ks_210}}
   \subfigure[]{%
   \includegraphics[width=5.5cm]{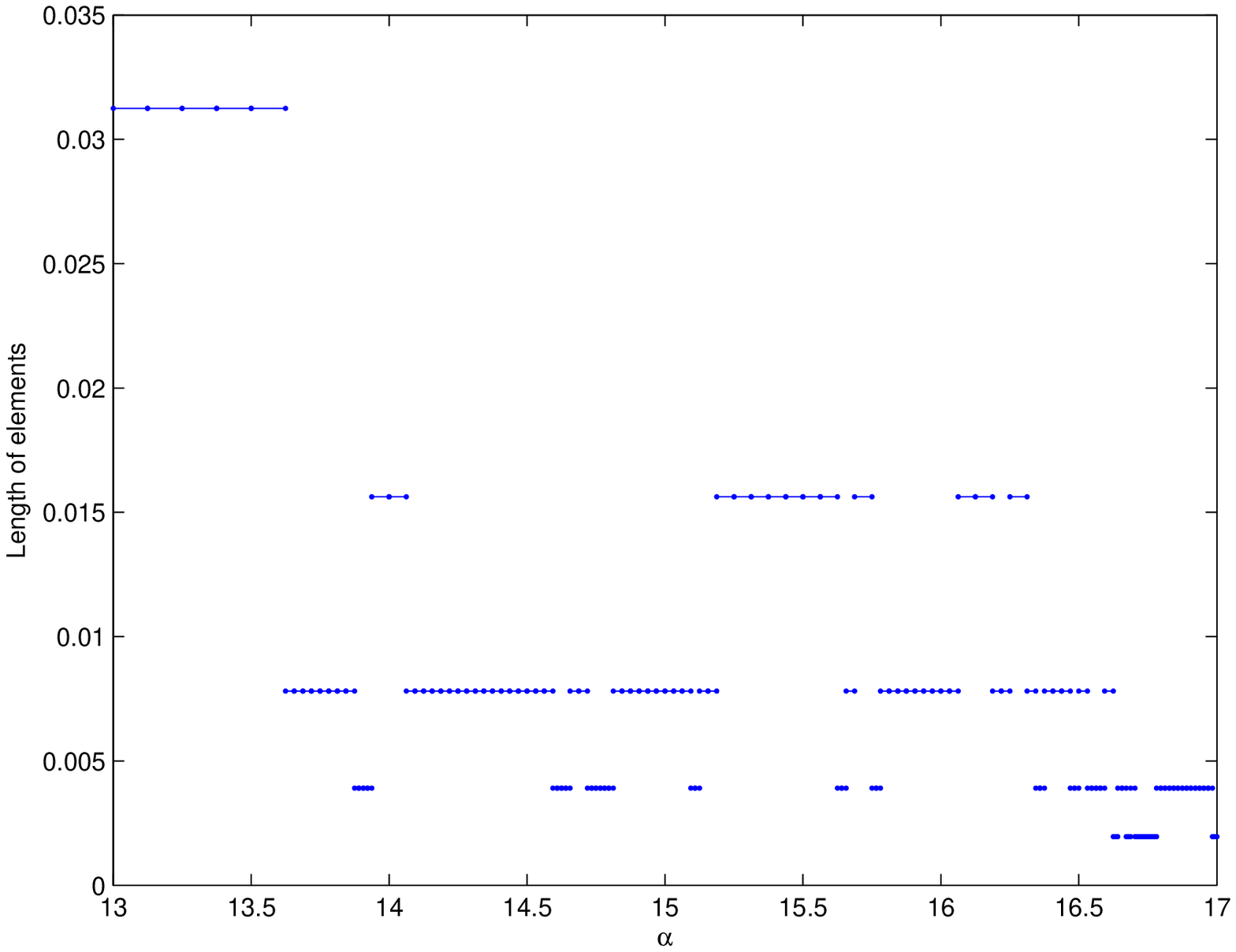}
   \label{fig:ks_end}}
\caption{Variance of solution to K-S equation at $x = 5.1296$ (left), Elements of  parameter space at $t = 10$ via AMR-CO with $d=11, TOL_1=0.1$(right)}
\label{fig:ks_1317}
\end{figure}

Figure \ref{fig:ex3_2} shows the evolution of the mean of K-S solution vis AMR-CO with $d=11, TOL_1 = 0.1$ compared with that of the reference. Figure \ref{fig:ex3_3} shows the comparison of the evolution of the variance of K-S solution between AMR-CO and the reference MC-SOBOL. We can see that although the behavior of the solution is very complicated, the AMR-CO still captures the major properties of the solution. The evolution variances at $x = 5.1296$ by different orders compared with the reference is shown in Figure \ref{fig:ks_210}. In this case, a more actuate model($d=11$) results in a better simulation. The distributions and length of the elements of the parameter space are presented in Figure \ref{fig:ks_end}. Since when $d=11$, there are $12$ collocation points in each element, in this case more smaller elements as well as more points are distributed in the bimodal solution regime. The final mesh has a lot of elements which demonstrate extreme sensitivity of the solution to the bifurcation parameter.

\subsection{Refinement in physical space}
\begin{figure}[htbp]
   \centering
   \subfigure[]{%
   \includegraphics[width=5.5cm]{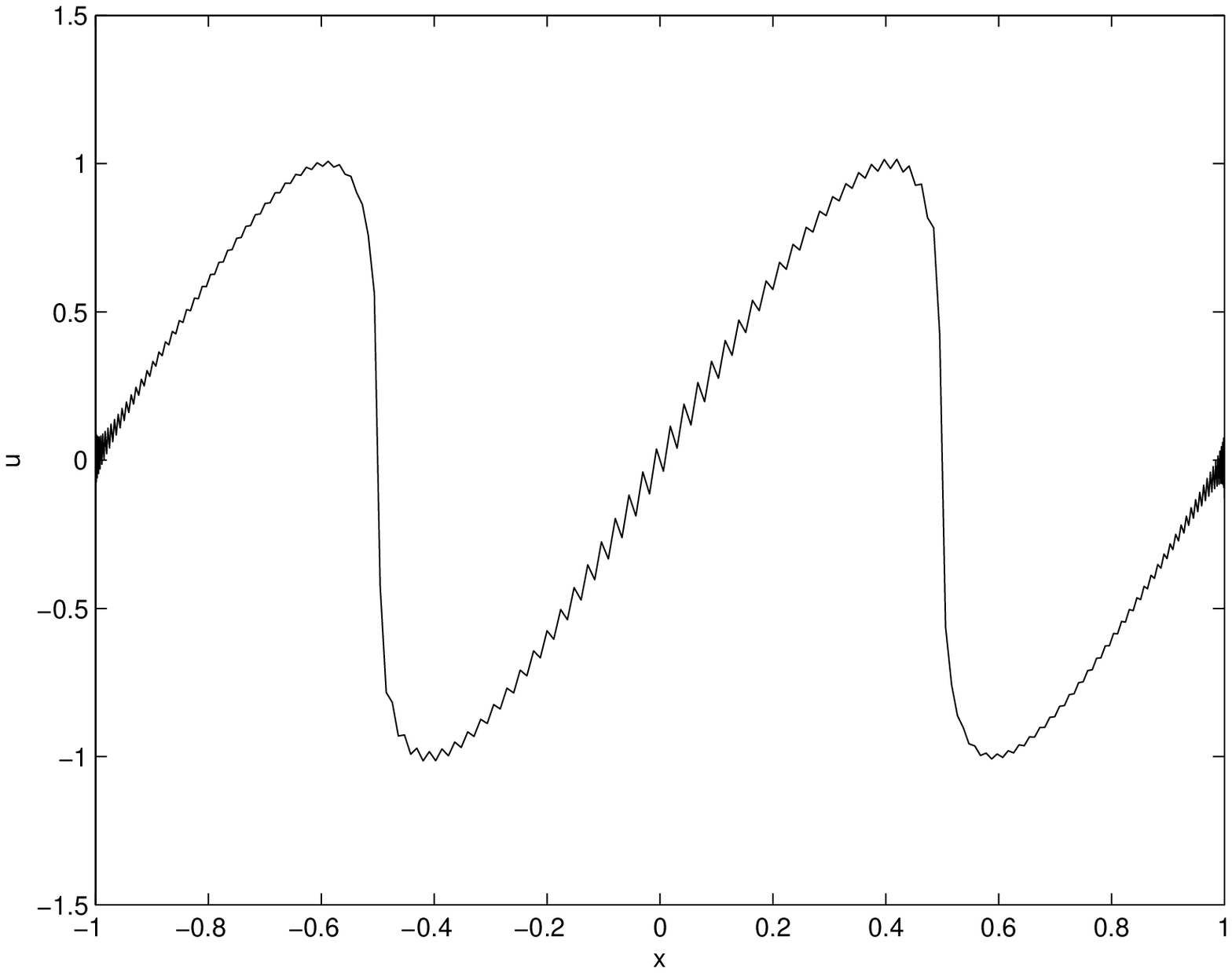}
   \label{fig:burgers_255}}
   \subfigure[]{%
   \includegraphics[width=5.5cm]{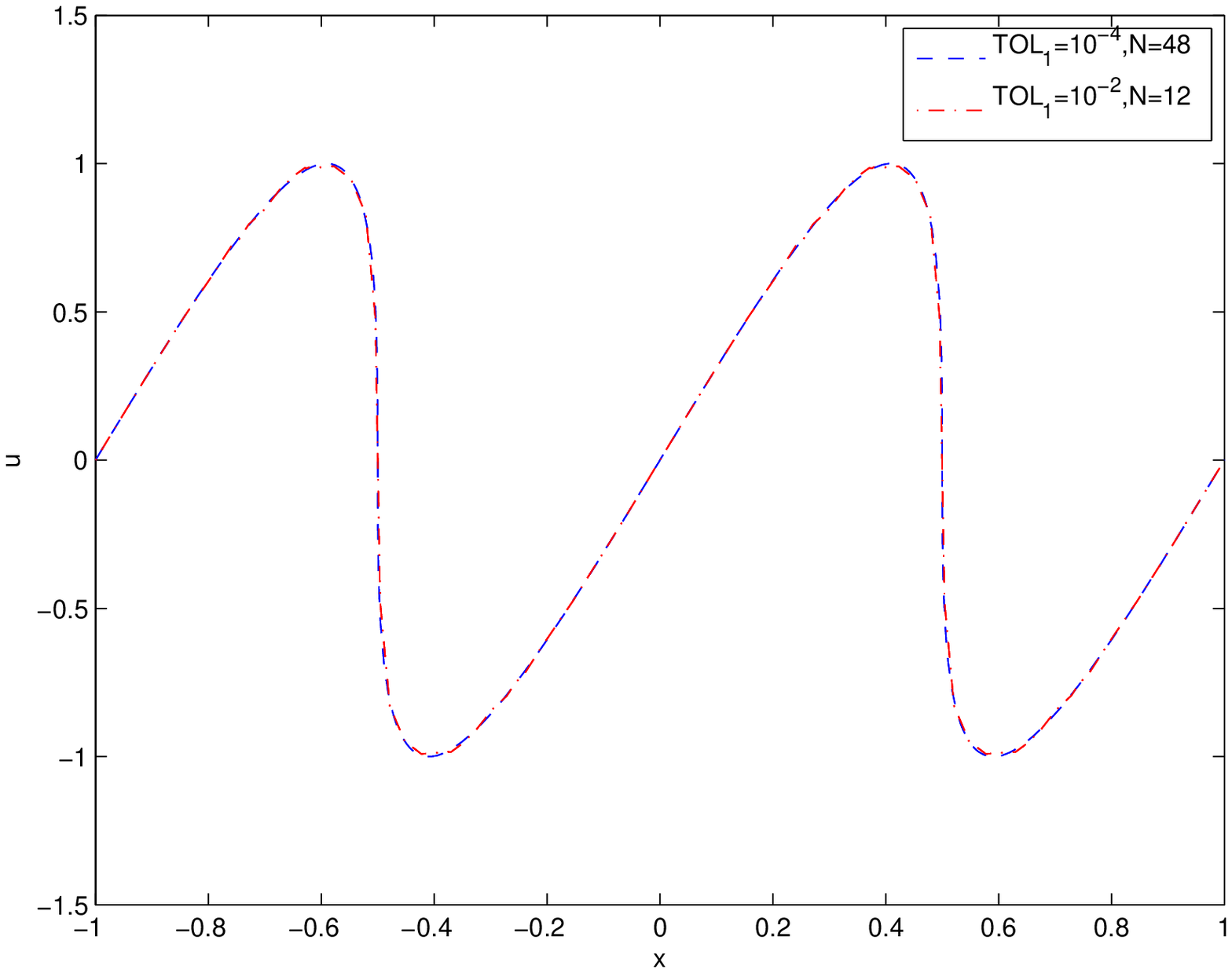}
   \label{fig:burgers_25_24}}
\caption{Solution of burgers equation by 256 Legendre collocation points(left) and mesh refinement algorithm with $p=5$, $TOL_1=10^{-2},10^{-4}$(right), here $\Delta t = 10^{-5}$}
\label{fig:burgers}
\end{figure}

\begin{figure}[htbp]
   \centering
   \includegraphics[width=12cm]{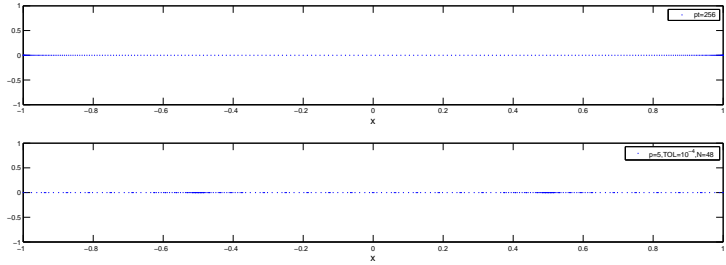}
   \caption{Distributions of 256 Legendre collocation points(top) and collocation point from mesh refinement algorithm with $p=5$, $TOL_1=10^{-2},10^{-4}$(bottom)}
\label{fig:burgers_dist_amr}
\end{figure}

Even though we have presented the formulation of our mesh refinement algorithm for random space, the same framework can accommodate also refinement in physical space. We will present here results for the location of shocks for the 1D Burgers equation (more elaborate examples will be presented in a forthcoming publication). The 1D inviscid Burgers equation reads 
 \begin{equation}\label{burgers}
u_t + u u_x = 0, \quad x\in(-1,1),
\end{equation}
with initial condition $u(x,0) = u_0(x)$ and boundary condition $u(-1,t) = u(1,t) = 0$. 
We use the initial condition:
\begin{equation}
u(x,0) = \sin(2\pi x), \quad x\in(-1,1).
\end{equation}
The solution develops two shocks in finite time. We employ a collocation method and solve the resulting system of ODEs by a fourth-order Runge-Kutta scheme with $\Delta t = 10^{-5}.$. Figure \ref{fig:burgers_255} shows the solution of \eqref{burgers} at $t=0.1592 (\approx 1 / \pi)$ using a global collocation method with
256 Gauss-Legendre collocation points. Because the distribution of the collocation points are concentrated near the boundaries of the domain as seen in top figure in Figure \ref{fig:burgers_dist_amr}, the solution shows oscillations in the areas where the gradient of the solution is large. Figure \ref{fig:burgers_25_24} shows the solution at $t=0.1592$ using the adaptive mesh refinement collocation algorithm. To resolve the solution we start with 8 elements and 6 collocation points in each element. As shown in Figure \ref{fig:burgers_dist_amr} the adaptive mesh refinement algorithm distributes more elements where the slope is large, which means more points near the location of the shocks. As a result the oscillations are eliminated. Note that the assignment of more points where a better resolution is needed stabilizes the numerical scheme. Even with large time step $\Delta t = 10^{-2}$, the mesh refinement algorithm is still stable, while direct collocation with 256 Gauss-Legendre point becomes unstable and oscillates a lot.

\section{Discussion and future work}\label{sec:conclusion}
We have presented a novel method for adaptive mesh refinement in the context of uncertainty quantification when a spectral decomposition of the solution is possible. The layering structure of the spectral representation of the solution can capture the transfer of activity across scales and thus can be utilized to detect when and where higher resolution is needed. This idea is inspired by \cite{Li2015} where the reduced model of the full system can be constructed. In this case the memory term of the reduced model captures the energy transfer between scales and provides a reliable indicator for the mesh refinement process. The proposed approach was implemented in the context of multi-element generalized polynomial chaos expansions for both Galerkin and collocation methods.

The same mesh refinement framework can be immediately extended to cases requiring refinement in physical space. However, unlike the random space where each element can be considered independently, when differentiating or integrating in physical space all elements are coupled due to the need of regularity of the solution across the elements. In the current work we have only included results from the application of the mesh refinement algorithm to the accurate location of shocks for the 1D Burgers equation. More elaborate examples (for more spatial dimensions) will be presented in a forthcoming publication.%

%
\section*{Acknowledgements} 
This material is based upon work supported by the U.S. Department of Energy Office of Science, Office of Advanced Scientific Computing Research, Applied Mathematics program, Collaboratory on Mathematics for Mesoscopic Modeling of Materials (CM4), under Award Number DE-SC0009280.
\bibliographystyle{elsarticle-num-names}
\bibliography{random}

\end{document}